\documentclass[a4paper,12pt]{article}
\usepackage{amssymb}
\usepackage{amsmath}
\usepackage{amscd}
\usepackage[arrow,matrix]{xy}
\usepackage[bookmarksnumbered,colorlinks,linkcolor=black,citecolor=black,urlcolor=black]{hyperref}
\usepackage{mathtools}
\usepackage{indentfirst} 

\usepackage{tikz-cd}
\usepackage{tikz}

\newtheorem{theo}{Theorem}[section]
\newtheorem{prop}[theo]{Proposition}
\newtheorem{lem}[theo]{Lemma}
\newtheorem{cor}[theo]{Corollary}
\newtheorem{defi}[theo]{Definition}
\newtheorem{rem}[theo]{Remark}
\newtheorem{exe}[theo]{Exercise}
\newtheorem{exa}[theo]{Example}

\newtheorem{conj}[theo]{Conjecture}
\newtheorem{ques}[theo]{Question}

\newcommand{\bthe}{\begin{theo}}
\newcommand{\ble}{\begin{lem}}
\newcommand{\bpr}{\begin{prop}}
\newcommand{\bco}{\begin{cor}}
\newcommand{\bde}{\begin{defi}}
\newcommand{\ethe}{\end{theo}}
\newcommand{\ele}{\end{lem}}
\newcommand{\epr}{\end{prop}}
\newcommand{\eco}{\end{cor}}
\newcommand{\ede}{\end{defi}}
\newcommand{\brem}{\begin{rem}}
\newcommand{\erem}{\end{rem}}

\newcommand{\bexe}{\begin{exe}}
\newcommand{\eexe}{\end{exe}}

\newcommand{\bexa}{\begin{exa}}
\newcommand{\eexa}{\end{exa}}

\newcommand{\bconj}{\begin{conj}}
\newcommand{\econj}{\end{conj}}
\newcommand{\bques}{\begin{ques}}
\newcommand{\eques}{\end{ques}}

\DeclareMathOperator{\pexp}{\mathnormal{p}-exp}

\def\fppf {{\rm fppf}}

\def \Et{{\rm Et}}

\def \Br {{\rm{Br}}}
\def \Brhat{\widehat{\mathrm{Br}}}

\def \si {{\sigma}}
\def \Ga {{\Gamma}}

\def \Pic {{\rm {Pic}}}

\def \Ker {{\rm{Ker}}}
\def \Coker {{\rm{Coker}}}
\def \Im {{\rm {Im}}}

\def \bH{{\bf H}}

\def \Spec {{\rm{Spec}}}
\def \dim {{\rm{dim}}}
\def \Hom {{\rm {Hom}}}
\def \End {{\rm {End}}}
\def \Pic {{\rm {Pic}}}

\def \cris {{\rm cris}}

\def\ov{\overline}

\def \Z {{\mathbb Z}}
\def \Q {{\mathbb Q}}
\def \F {{\mathbb F}}

\def \red {{\rm red}}

\def \rk {{\rm{rk}}}

\def \H {{\rm H}}
\def\G{{\mathbb G}}

\def\sZ{{\cal Z}}
\def\sY{{\cal Y}}

\def\lra{\longrightarrow}

\def\NS{{\rm NS\,}}
\def\skew{{\rm skew}}

\def\O{{\cal O}}

\def\sX{{\cal X}}
\def\sY{{\cal Y}}

\def\si{\sigma}

\def\Ga{\Gamma}

\def\et{{\rm{\acute et}}}

\def\perf{{\rm perf}}

\def\bL{{\bf L}}
\def\bU{{\bf U}}
\def\bG{{\bf G}}
\def\bW{{\bf W}}

\def\E{{\mathbb E}}

\def\Nyg{{{\rm F}^1_N{}}}

\newcommand{\W}{\mathbb{W}}

\definecolor{blu}{RGB}{0,114,142}
\definecolor{rosso}{RGB}{183,42,58}

\usepackage[textheight=220mm,textwidth=150mm]{geometry}

\title{Brauer groups of abelian varieties over fields of finite characteristic}
\author{Livia Grammatica, Alexei N. Skorobogatov, and Yuan Yang}

\date{\today}
\begin{document} 
\maketitle

\begin{abstract}
We study the Brauer group of an abelian variety $A$ over an algebraically closed field of
characteristic $p$ focusing on the $p$-primary torsion, the key part of which is a certain 
quasi-algebraic unipotent group $\bU_A$. We determine its dimension 
and obtain a sharp upper bound for its $p$-exponent. The isogeny class of
$\bU_A$ is classified for abelian varieties $A$ of dimension at most 3.
For principally polarised abelian varieties we compute the dimension of the $p$-torsion subgroup
of $\bU_A$ in terms of the Ekedahl--Oort type of $A$.  
\end{abstract}

\tableofcontents

\section*{Introduction}

Let $k$ be an algebraically closed field of characteristic $p\geqslant 0$. 
In this paper we are interested in the following question:
what is the Brauer group of an abelian variety $A$ over $k$? In particular, 
what is the $p$-primary torsion subgroup $\Br(A)[p^\infty]$ when $p>0$?

Let us start by recalling the structure of the Brauer group of an arbitrary smooth proper variety $X$
over $k$.
In \cite[III]{GB} Grothendieck proved that $\Br(X)$ is a torsion group
and computed $\Br(X)[\ell^\infty]$
for primes $\ell\neq p$. The exactness of the Kummer sequence with 
coefficients $\mu_{\ell^n}$ for the \'etale topology gives an isomorphism
$$\Br(X)[\ell^\infty]\cong (\Q_\ell/\Z_\ell)^{b_2-\rho}\oplus \H^3_\et(X,\Z_\ell(1))[\ell^\infty],$$ where $b_2=\dim_{\Q_\ell}\H^2_\et(X,\Q_\ell)$ is the second $\ell$-adic Betti number,
$\rho=\rk\,\NS(X)$ is the Picard number, and 
$\H^3_\et(X,\Z_\ell(1))[\ell^\infty]$ is finite, see \cite[Prop.~5.2.9]{CTS21}. 
When $X=A$ is an abelian variety of dimension $g>0$, 
we have $b_2=g(2g-1)$ and $\H^3_\et(A,\Z_\ell(1))[\ell^\infty]=0$.
In particular, the abelian group $\Br(A)[\ell^\infty]$ is an invariant of the isogeny class of $A$.

The $p$-primary torsion subgroup $\Br(X)[p^\infty]$ in characteristic $p>0$
has more interesting structure. 
Let $W=W(k)$ be the ring of Witt vectors of $k$. 
Define $\H^i(X,\Z_p(1))=\varprojlim\H^i_\fppf(X,\mu_{p^n})$. 
By Illusie, we have $\H^i(X,\Z_p(1))/{\rm tors}\cong\Z_p^{r_i}$,
where $r_i$ is the multiplicity of slope 1 
of the crystalline cohomology group $\H^i_\cris(X/W)$, see \cite[Thm.~II.5.5 (5.5.3)]{Ill}. 
By Illusie--Raynaud \cite[Thm.~IV.3.3 (b), Cor.~IV.3.6]{Ill-Ray}
it is known that $\H^i(X,\Z_p(1))[p^\infty]$ is annihilated by a power of $p$.
The exactness of the Kummer sequence with 
coefficients $\mu_{p^n}$ for the fppf topology gives an isomorphism
\begin{equation}
\Br(X)[p^\infty]\cong (\Q_p/\Z_p)^{r_2-\rho}\oplus \H^3(X,\Z_p(1))[p^\infty], \label{sasha}
\end{equation}
see \cite[Theorem A.1]{P}. In contrast with the $\ell$-adic case, the group 
$\H^3(X,\Z_p(1))[p^\infty]$ can be infinite.
A phenomenon specific to finite characteristic is the existence of a natural commutative 
quasi-algebraic group over $k$ whose group of $k$-points is $\H^3(X,\Z_p(1))[p^\infty]$.
Its connected component, which is unipotent, will be denoted by $\bU_X$.
This goes back to Milne \cite[Cor.~2.7 (a)]{Milne}, see \cite[Cor.~2.10]{Ber81},
\cite[Lemme IV.3.2.2]{Ill-Ray}, \cite[Lemma 1.8]{Milne2}.
Quasi-algebraic groups were introduced by Serre \cite[\S 1.2]{Serre60} as equivalence classes
of algebraic groups with a common set of $k$-points, where two groups $G$ and $H$ are called equivalent if there is a morphism $G\to H^{(p^n)}$ which is a bijection on $k$-points.

When $X=A$ is an abelian variety, the attached quasi-algebraic group is connected and so is equal
to $\bU_A$, see
\cite[Thm.~1.5]{Yang} based on  \cite[Thm.~II.3.4]{Ill-Ray} (``survie du c{\oe}ur"). 
In particular, the group of $k$-points $\bU_A(k)\cong\H^3(A,\Z_p(1))[p^\infty]$, the dimension
$\dim(\bU_A)$, and the $p$-exponent $\pexp(\bU_A)$ 
(the least integer $m$ such that $p^m\bU_A=0$)
are all well-defined. By general theory \cite[VII.10]{Serre},
$\bU_A$ is isogenous to a product of truncated Witt (quasi-algebraic) groups 
$\prod_{i=1}^r\bW_{n_i}$, 
for some positive integers $n_1\leqslant \ldots \leqslant n_r$.
We have 
$$\dim(\bU_A)=\sum_{i=1}^r n_i, \quad \dim(\bU_A[p])=r, \quad\pexp(\bU_A)=n_r.$$
Note that if $\pexp(\bU_A)=1$, then $\bU_A\simeq \bG_{a,k}^d$, see 
\cite[VII.2, Prop.~11]{Serre}. 

The guiding question of this paper is: what sequences $(n_1,\ldots,n_r)$ occur for abelian varieties
of dimension $g$ over an algebraically closed field $k$ of positive characteristic $p$?

Using results of Crew and deep results of Ekedahl,
together with the fact that the groups $\H^i_\cris(A/W)$ are torsion-free
and the Hodge-to-de Rham spectral sequence for $A$ degenerates at the first page, one obtains
\cite[Theorem 1.1, 1.6]{Yang} 
\begin{equation}
\dim(\bU_A)=\frac{g(g-1)}{2}-\sum_{0\leqslant \lambda<1}(1-\lambda)m_\lambda, \label{dim}
\end{equation}
where $g=\dim(A)$ 
and $m_\lambda$ is the multiplicity of slope $\lambda$ in $\H^2_\cris(A/W)$. 
In particular, $\dim(\bU_A)$ depends only on the slopes of $A$
and thus is an invariant of the isogeny class of $A$. 
The meaning of the terms of (\ref{dim}) is the following: 
$\dim_k\H^0(A,\Omega^2_{A/k})=g(g-1)/2$ 
is the dimension of the tangent space of the formal Brauer group $\widehat{\Br}(A)$,
the sum is the dimension of the maximal quotient of $\widehat{\Br}(A)$ of finite height,
and $\dim(\bU_A)$ is the dimension of the unipotent part of $\widehat{\Br}(A)$,
see Section \ref{formal}. Thus (\ref{dim})
implies the smoothness of the formal Brauer group of any abelian variety over $k$
(Corollary \ref{cor3.3}). Note that
smoothness of $\widehat{\Br}(A)$ is a particular case of a general criterion of smoothness
of Artin--Mazur formal groups obtained by the first named author in \cite{Gr2}.

From (\ref{dim}) we see that $\dim(\bU_A)\leqslant g(g-1)/2$ with equality if and only if $g=1$
or $A$ is supersingular (Remark \ref{ferry}).
For $g=2$, it gives $\bU_A=0$ unless $A$ is supersingular, in which case
$\bU_A\cong\bG_{a,k}$, a fact that goes back to Illusie  
whose computation was based on the de Rham--Witt cohomology \cite[II.7]{Ill}. In general,
it is known that $\bU_A=0$ when $\H^2(A,W{\mathcal O}_A)$ is a finitely generated $W$-module, 
for example, when $A$ is ordinary (i.e., the $p$-rank is $g$) or almost ordinary
(i.e., the $p$-rank is $g-1$), see \cite[Corollary 6.3.16]{Ill83}, \cite[Theorem 3.7]{Yang}.
In fact, $\bU_A=0$ if and only if $A$ is ordinary or almost ordinary, see Remark \ref{ferry}.

Let us now describe the main results of this paper and give its brief outline.

\begin{itemize}

\item
Developing the approach of \cite{Gr}
based on the Nygaard filtration, we express $\bU_A$ in terms of the F-crystal 
$\H^2_\cris(A/W)$, see Corollary \ref{xx1}.
This gives a quick proof that $\dim(\bU_A)$ depends only on the isogeny class of $A$
(Proposition \ref{iso-inv}).
From this we obtain a proof of formula (\ref{dim}) that uses only the classical Dieudonn\'e theory and
so avoids the deep theory of Ekedahl, see Theorem~\ref{3.1}. At the end of \S \ref{formal}
we give a different, short proof of (\ref{dim}) based on the smoothness of the formal Brauer group
of $A$.

\item
We interpret $\Br(A)$ in terms of homomorphisms of truncated
Dieudonné modules. Namely, for $p\neq 2$ there is a natural isomorphism of 
$\bU_A[p^n]$ and the quasi-algebraic group
$$\Hom_\E(M_A^\vee/p^n,M_A/p^n)^\skew/\left(\Hom_\E(M_A^\vee,M_A)^\skew/p^n\right),$$
which is a subquotient of the quasi-algebraic group $\Hom_{W(k)}(M_A^\vee/p^n,$ $M_A/p^n)$, see Corollary \ref{2.14}. Here $M_A$ is the Dieudonné module of $A$, 
$\mathbb{E}$ is the Dieudonné ring (see \S \ref{Not}) and 
`$\skew$' stands for the subgroup of anti-self-dual elements.
This gives a different way to reduce the determination of $\bU_A$ to
a computation with the Dieudonn\'e module $M_A$. 

\item
We prove that the $p$-exponent of $\bU_A$ is at most $g-1$, where $g=\dim(A)$ and $p\neq 2$,
see Theorem \ref{ss}. We show that the bound is sharp by exhibiting,
for every $g\geqslant 1$, a supergeneral (i.e.~supersingular with $a$-number 1)
abelian variety $A$
of dimension $g$ such that $\bU_A$ is isogenous to $\prod_{n=1}^{g-1}\bW_n$. Combining this with
a result of D'Addezio we obtain that if $A$ is an abelian variety over a field $k_0$
 finitely generated over $\F_p$, and $k=\ov{k_0}$, then the transcendental Brauer group
$\Im[\Br(A)\to\Br(A_k)]$ is a direct sum of a finite group and an abelian group annihilated by
$p^{g-1}$ (Corollary \ref{f.g.}). 
Moreover, if $A$ is supersingular, then the transcendental Brauer group of $A$ is annihilated by
$p^{g-1}$, whereas if $A$ is ordinary or almost ordinary, then this group is finite.

\item
We determine the isogeny class of $\bU_A$ for arbitrary abelian threefolds (assuming $p\neq 2$): $\bU_A\simeq\bG_{a,k}^n$, where $0\leqslant n\leqslant 3$,
unless $A$ is supergeneral, in which case
$U_A$ is isogenous to $\bG_{a,k}\times \bW_2$, see Proposition \ref{NPPabe3table}.
In particular, the isogeny class of $\bU_A$ is {\em not} an invariant of the isogeny class of $A$.

\item
We compute $\dim(\bU_A[p])$ for principally polarised abelian varieties $A$ 
of arbitrary dimension $g\geqslant 1$ in terms of the Ekedahl--Oort type of $A[p]$, 
see Theorem \ref{E-O}. From this we deduce that in the principally polarised case we have
$\dim(\bU_A)\leqslant r(r+1)/2$ 
and $\pexp(\bU_A)\leqslant r$, where $r=\dim(\bU_A[p])$, see Corollary \ref{d1}.
\end{itemize}

The work on this paper started when A.S. was a Shapiro Visitor
in Penn State. He is very grateful to Yuri Zarhin 
 for many helpful discussions without which this paper would have never been written.
L.G. received funding from the EU Horizon Europe research
and innovation programme under the Marie Skłodowska-Curie grant agreement 101126554.
The work of Y.Y. was supported by the London School of Geometry and Number Theory
through the Engineering and Physical Sciences Research Council 
grant EP/S021590/1.  The authors are grateful to the Bernoulli Center 
at EPFL, where a large part of this work was carried out,
for excellent working conditions and support.
We would like to thank James Borger, George Boxer, Daniel Bragg, Marco D'Addezio, Vladimir Drinfeld, Mikhail Kapranov, and Alexander Petrov for their questions and helpful discussions.

\section{Preliminaries}
\subsection{Notation} \label{Not}

In this paper $k$ is an algebraically closed field of characteristic $p>0$.
We write $W=W(k)$ for the ring of Witt vectors, and 
denote by $\si\colon W\to W$ the canonical lifting of the Frobenius $x\to x^p$ on $k$.

All group schemes in this paper are commutative. 
For a finite flat commutative group $k$-scheme $G$ we denote by $G^\vee$ its Cartier dual.
If $A$ is an abelian variety, we write $A^\vee$ for the dual abelian variety. 

The Dieudonn\'e ring $\E=W_\si[F,V]$ is the ring of non-commutative polynomials in $F$ and $V$ with coefficients in $W$ subject to the relations
$FV=VF=p$, $F\lambda=\lambda^\si F$, $V\lambda=\lambda^{\si^{-1}}V$, where $\lambda\in W$.
An $\E$-module $M$, which is free and finitely generated as a $W$-module, will be called
a Dieudonn\'e module. The dual Dieudonn\'e module is defined as $M^\vee=\Hom_W(M,W)$
with the action of $\E$ on $f\in M^\vee$ given by the rules $Ff(x)=\si f(Vx)$ and $Vf(x)=\si^{-1}f(Fx)$, 
for any $x\in M$.

\subsection{Quasi-algebraic and pro-algebraic groups}

Quasi-algebraic groups were introduced by Serre \cite{Serre60} as a way to have an abelian category of group schemes where infinitesimal groups and radicial morphisms are ignored. This is useful whenever one is interested in working with $k$-points of algebraic groups, rather than algebraic groups themselves, as in our case. Serre's construction can be interpreted as restricting a group scheme over $k$ to the category of perfect schemes over $k$.

Let $\pi\colon X\to\Spec(k)$ be a proper and smooth
scheme over $k$. Let $\Spec(k)_\perf$ be the {\em perfect} site of $k$. 
The objects of $\Spec(k)_\perf$
are perfect $k$-schemes (that is, $k$-schemes such that the Frobenius map is an isomorphism)
equipped with the \'etale topology. Let $(X/k)_\perf$ be the 
site such that the objects are pairs $(T,Y)$, where $T$ is a perfect
$k$-scheme and $Y$ is \'etale over $X\times_k T$, equipped with the \'etale topology. 
By \cite[2.7]{Ber81} we have a commutative diagram of topoi, where $X_\Et$ is the big \'etale site of $X$,
and $\alpha$ and $\beta$ are natural restriction maps:
$$\xymatrix{X_\fppf\ar[r]^{\beta_X}\ar[d]_{\pi_\fppf}&X_\Et\ar[r]^{\alpha_X \ \ \ }\ar[d]_{\pi_\Et}
&(X/k)_\perf\ar[d]_\pi\\
\Spec(k)_\fppf\ar[r]^{\beta_k}&\Spec(k)_\Et\ar[r]^{\alpha_k \ \ } &\Spec(k)_\perf}$$
By definition, perfect $k$-groups are the group objects in $\Spec(k)_\perf$. 
The {\em perfection}
$G^\perf$ of a group $k$-scheme $G$ is defined as the inverse limit 
$\varprojlim G_\red^{(p^{-n})}$ of iterated Frobenius twists
along the Frobenius morphisms $F\colon G_\red^{(p^{-(n+1)})}\to G_\red^{(p^{-n})}$
for $n\geqslant 0$.
The restriction $\alpha_k\circ\beta_k$ of the sheaf on $\Spec(k)_\fppf$ represented by $G$ to 
$\Spec(k)_\perf$ is a sheaf represented by $G^\perf$. The functors $\alpha_{k,*}$ and 
$\alpha_{X,*}$ are exact ({\em loc.~cit.}).

Let ${\mathcal G}[p^\infty]$ be the category of perfect commutative affine $k$-groups
annihilated by some power of $p$ that are perfections of group $k$-schemes of finite type.
Let ${\mathcal P}$ be the category of pro-objects of ${\mathcal G}[p^\infty]$.
The objects of ${\mathcal G}[p^\infty]$ (respectively, ${\mathcal P}$)
can be described as {\em quasi-algebraic} (respectively, {\em pro-algebraic}) $k$-groups
in the sense of Serre \cite[\S 1, \S 2]{Serre60}. Serre shows that ${\mathcal G}[p^\infty]$
and ${\mathcal P}$ are abelian categories, where subobjects are closed subgroups
\cite[\S 1.2, Prop.~5, \S 2.4, Prop.~7]{Serre60},
see also \cite[2.2 (a)]{Ber81}. An example of an object of ${\mathcal G}[p^\infty]$
is ${\bf W}_n=(\W_n)^\perf$, where $\W_n$ is the group $k$-scheme of Witt vectors of length $n$.
In particular, we have $\bG_{a,k}:={\bf W}_1=(\G_{a,k})^\perf$.
An example of an object of 
${\mathcal P}$ is ${\bf W}=\W^\perf$, see \cite[p.~13, Exemple 3]{Serre60}.

\subsubsection*{Unipotent quasi-algebraic groups from F-crystals}

Let $L$ be an F-crystal, i.e.~a free finitely generated $W$-module with an injective
 $\si$-linear endomorphism $F$. Define 
$$\Nyg(L)=\{x\in L\text{ such that }p|F(x)\}.$$
In particular, we have $pL\subset \Nyg(L)\subset L$.

The F-crystal $L$ has a natural structure of a pro-algebraic group $\bL$ such that $\bL(k)=L$,
so that $\bL$ is isomorphic to ${\bf W}^n$ where $n=\rk_W(L)$.
Then $F$ naturally gives rise to a $k$-endomorphism $\bL\to\bL$
also denoted by $F$. It is the relative Frobenius 
$F_{\bL/k}\colon \bL\to\bL$ followed by a linear map $\bL\to\bL$.

Thus $\Nyg(\bL)$ is a closed subgroup of $\bL$, and $F/p-1\colon \Nyg(\bL)\to \bL$ 
is a morphism in the category of pro-algebraic $k$-groups ${\mathcal P}$. 
But ${\mathcal P}$ is an abelian category, so 
\begin{equation}
\bU_L:=\Coker[F/p-1\colon \Nyg(\bL)\to \bL] \label{def1}
\end{equation}
is a pro-algebraic group.

\bpr \label{s1}
Let $L$ be an F-crystal. Then $\bU_L$ is a connected unipotent quasi-algebraic group.
\epr
{\em Proof.} Let $n$ be the rank of the free finitely generated $W$-module $L$.
By a standard argument \cite[Lemme II.5.6]{Ill}, the cokernel of $F/p-1\colon \Nyg(L)\to L$
is a torsion abelian group.
We claim that it is annihilated by $p^r$ for some $r\geqslant 1$. From this it follows that 
$\bU_L$ is the quotient of the connected unipotent quasi-algebraic group $k$-scheme ${\bf W}_r^n$ 
by a closed subgroup. This implies the proposition.

It remains to prove the claim. We proceed as in \cite[Lemma 5.1.3]{Gr}. Since $L/\Nyg(L)$ is killed by $p$, it is enough to prove
that the cokernel of $F-p\colon L\to L$ has finite $p$-exponent. By the Manin--Dieudonn\'e
classification, it is enough to prove this for an F-crystal $M$ with a $W$-basis $e_1,\ldots,e_r$
such that $Fe_i=e_{i+1}$ for $i=1,\ldots,r-1$ and $Fe_r=p^se_1$, where $r>0$ and $s\geqslant 0$
are coprime integers. There are three possible cases:

(1) If $r=s=1$, then $F-p=p(F/p-1)$, and $(F/p-1)(\lambda e_1)=(\si(\lambda)-\lambda)e_1$ is surjective by \cite[Lemme II.5.3]{Ill}.

(2) If $s>r$, then $p^rM\subset(F-p)M$. Indeed, if $a\in p^rM$, then the infinite sum
$-\sum_{i\geqslant 1}F^{i-1}p^{-i}a$ converges to some $x\in M$, and we have $a=Fx-px$.

(3) If $s<r$, then $p^sM\subset(F-p)M$. Indeed, if $a\in p^sM$, then $\sum_{i\geqslant 0}F^{-i}p^ia$
converges to some $x\in M$, and we have $a=Fx-px$.  \hfill $\Box$

\bde \label{defU_L}
We call  $\bU_L$ defined in $(\ref{def1})$
the connected unipotent quasi-algebraic group attached to 
the F-crystal $L$.
\ede

The $p$-exponent of a unipotent quasi-algebraic group $\bU$ is the least
non-negative integer $n$ such that $p^n\bU=0$.

\bco \label{s3}
Define $\bH_{L,n}=\Ker[F/p-1\colon\Nyg(\bL)/p^n\to \bL/p^n]$. Then we have an exact sequence
of quasi-algebraic $k$-groups
\begin{equation}
0\to \bL^{F=p}/p^n \to \bH_{L,n}\to \bU_L[p^n]\to 0. \label{s2}
\end{equation}
In particular, for $n$ greater or equal to the $p$-exponent of $\bU_L$ the quasi-algebraic group
$\bU_L$ is naturally isomorphic to the cokernel of $\bL^{F=p}/p^n \to \bH_{L,n}$.
\eco
{\em Proof.} We have an exact sequence of pro-algebraic groups
$$0\to \bL^{F=p}\to  \Nyg(\bL)\xrightarrow{F/p-1} \bL\to \bU_L\to 0.$$
Breaking it into two short exact sequences and applying the snake lemma
to the diagrams formed by multiplication by $p^n$, we obtain (\ref{s2}). \hfill $\Box$

\medskip

Note that $\bL^{F=p}$ is finitely generated as a $\Z_p$-module \cite[Lemme II.5.11]{Ill},
so $\bL^{F=p}/p^n$ is finite. Thus the map $\bH_{L,n}\to \bU_L[p^n]$ is an isogeny
of quasi-algebraic groups.

\subsection{Quasi-algebraic groups from flat and crystalline cohomology}

The sheaves $W_n\Omega^i$ on $X_\et$ can be lifted to sheaves 
on $(X/k)_\perf$, again denoted $W_n\Omega^i$, see \cite[p.~192]{Ill-Ray}.
Formation of $W_n\Omega^i$ commutes with perfect base change \cite[Prop.~I.1.9.2]{Ill}.
From this one deduces that the sheaves $R^j\pi_*(W_n\Omega^i)$ 
on $\Spec(k)_\perf$ belong to
${\mathcal G}[p^n]$, that is, are represented by quasi-algebraic groups,
see \cite[Lemme IV.3.2.2]{Ill-Ray}. Using the notation of \cite{Ill-Ray}
we denote these quasi-algebraic groups by
$\bH^j(X,W_n\Omega^i)$. 
Considering the complex $W_n\Omega_{X/k}^\bullet$
as an object of the derived category of sheaves on $(X/k)_{\perf}$, we
denote the quasi-algebraic group representing
$R^j\pi_{*}(W_n\Omega_{X/k}^\bullet)$ by $\bH^j_\cris(X/W_n)$.
Passing to the limit we obtain pro-algebraic $k$-groups 
$\bH^j(X,W\Omega^i)$ and  $\bH^j_\cris(X/W)$.
Since $k$ is algebraically closed, the group of $k$-points of
$\bH^j(X,W\Omega^i)$ (respectively, $\bH^j_\cris(X/W)$) is
$\H^j(X,W\Omega^i)$ (respectively, 
$\H^j_\cris(X/W)$), and the same holds with $W$ replaced by $W_n$.  
Moreover, for a perfect $k$-algebra $k'$ we have 
$$\bH_\cris^{i}(X/W)(k')=\Ga(\Spec(k'),R^j\pi_{*}(W\Omega_{X/k}^\bullet))= \H_\cris^{i}(X/W)\otimes_W W(k').$$
When $\H_\cris^{i}(X/W)$ is torsion-free, this implies that $\bH_\cris^{i}(X/W)$ 
is the pro-algebraic group ${\bf W}^n$ for some $n\geqslant 0$. 

On all of the above pro-algebraic groups, the Frobenius $F$ acts as an endomorphism over $k$.

\medskip

A fundamental result of Milne 
\cite[Lemma 1.8]{Milne2} says that the presheaf $T\mapsto \H^j_\fppf(X_T,\mu_{p^n})$ on $\Spec(k)_\perf$ is represented by a quasi-algebraic group, which we denote by $\bH^j(X,\mu_{p^n})$. 
Define $\bH^j(X,\Z_p(1))=\varprojlim \bH^j(X,\mu_{p^n})$. 
As $k$ is algebraically closed, the group of $k$-points of the pro-algebraic group
$\bH^j(X,\Z_p(1))$ is $\H^j(X,\Z_p(1)) :=\varprojlim \H^j_\fppf(X,\mu_{p^n})$.
By work of Illusie and Illusie--Raynaud \cite[Thm.~IV.3.3]{Ill-Ray} there are canonical isomorphisms
$$\bH^j(X,\Z_p(1))\cong\bH^{j-1}(X,W\Omega^1_{\rm log})\cong \bH^{j-1}(X,W\Omega^1)^{F=1}.$$

\medskip

As in \cite[Def.~8.1]{BMS} we define the first term $\Nyg W\Omega^\bullet$
of the Nygaard filtration on the de Rham--Witt complex $W\Omega^\bullet$ as the complex
$$0\to VWO_X\to W\Omega^1_X\to W\Omega^2_X\to\ldots$$
Thus there is an exact sequence of complexes of pro-\'etale sheaves
\begin{equation}
0\to \Nyg W\Omega^\bullet_X\to W\Omega_X^\bullet\to O_X \to 0, \label{n3}
\end{equation}
hence a long exact sequence of pro-algebraic groups
\begin{equation}
\ldots\to \bH^{j-1}(X,O_X)\to \Nyg\bH^j_\cris(X/W)\to \bH^j_\cris(X/W)\to \bH^{j}(X,O_X)\to\ldots \label{n1}
\end{equation}
One can define maps $F/p:\Nyg{}\bH_\cris^i(X/W)\to\bH_\cris^i(X/W)$ making the following diagram commute 
$$\xymatrix{\Nyg{}\bH_\cris^i(X/W)\ar[r]^{F/p}\ar[d]^1&\bH_\cris^i(X/W)\ar[d]^p\\
\bH_\cris^i(X/W)\ar[r]^F&\bH_\cris^i(X/W)}$$
By \cite[Prop.~8.4]{BMS} the map $F/p-1$ fits into an exact sequence
\begin{equation}
0\to W\Omega^1_{X,\log}[-1]\to \Nyg W\Omega^\bullet_X\xrightarrow{F/p-1}
W\Omega_X^\bullet\to 0 \label{n5}
\end{equation} 
This gives a long exact sequence of pro-algebraic groups
\begin{equation}\ldots\to \bH^j(X,\Z_p(1))\to \Nyg\bH^j_\cris(X/W)\xrightarrow{F/p-1}
\bH^j_\cris(X/W)\to \bH^{j+1}(X,\Z_p(1))\to\ldots \label{n2}
\end{equation}
see also \cite[Thm.~7.3.5]{bhatt-lurie} and \cite[Thm.~4.5]{Gr} for a different approach.

\bde
A smooth proper $k$-variety $X$ over an algebraically closed field $k$ of characteristic $p$ 
is called {\bf straight} or {\bf Mazur--Ogus}
if $\bH^i_\cris(X/W)$ is torsion-free for all $i\geqslant 0$, and the Hodge-to-de Rham spectral sequence of 
$X$ degenerates at the first page.
\ede

For example, abelian varieties and K3 surfaces are straight. 
For some properties of straight varieties, see \cite[\S 1]{Ogus}. 

\bpr \label{fppfsplit}
A straight variety $X$ over an algebraically closed field $k$ of characteristic $p$ has the following properties.

\smallskip

{\rm (i)} There is an exact sequence of pro-algebraic groups
$$0\to\Nyg{}\bH_\cris^{i}(X/W)\to \bH_\cris^{i}(X/W)\to \bH^i(X,\O_X)\to 0.$$ Here
$\Nyg{}\bH_\cris^{i}(X/W)\subset \bH_\cris^{i}(X/W)$ is the closed subgroup given by the condition that $p$ divides $F(x)$.
The map $$F/p-1\colon \Nyg{}\bH_\cris^{i}(X/W)\to \bH_\cris^{i}(X/W)$$ is the obvious map. 

\smallskip

{\rm (ii)} There is an exact sequence of pro-algebraic groups
$$
0\to\bH^i(X,\Z_p(1))/\mathrm{tors}\to\Nyg{}\bH_\cris^{i}(X/W)
\xrightarrow{\frac{F}{p}-1}\bH_\cris^{i}(X/W)\to\bH^{i+1}(X,\Z_p(1))[p^\infty]\to0.
$$
In particular, $\bH^{i+1}(X,\Z_p(1))[p^\infty]$ is the connected unipotent quasi-algebraic group
attached to the $F$-crystal $\bH_\cris^{i}(X/W)$.

{\rm (iii)} There is a long exact sequence of quasi-algebraic groups
$$\ldots\to\bH^{i}(X,\mu_{p^n})\to \Nyg{}\bH_\cris^{i}(X/W)/p^n
\xrightarrow{\frac{F}{p}-1}\bH_\cris^{i}(X/W)/p^n\to\bH^{i+1}(X,\mu_{p^n})\to\ldots$$
\epr
{\em Proof.} (i) The map $\bH_\cris^{i}(X/W)\to \bH^i(X,\O_X)$ is the composition of 
the natural maps
$$\bH_\cris^{i}(X/W)\xrightarrow{\alpha}\bH_{\rm dR}^{i}(X/W)\xrightarrow{\beta} 
\bH^{i}(X,\O_X).$$
Here $\alpha$ is surjective since $\bH_\cris^{i+1}(X/W)$ is torsion-free, and $\beta$
is surjective by the degeneration of the Hodge-to-de Rham spectral sequence.
Thus (\ref{n1}) gives the exact sequence in (i). The second property is a consequence
of the Mazur--Ogus theorem, see also \cite[Prop.~5.2.2]{Gr}.

(ii) The exact sequence 
follows from (\ref{n2}) since $\Nyg{}\bH_\cris^{i}(X/W)$ and $\bH_\cris^{i}(X/W)$ 
are both torsion-free, while
the cokernel of $F/p-1$ is torsion, by a standard argument in the theory of F-isocrystals 
\cite[Lemme II.5.6]{Ill}. Since $\bH_\cris^{i}(X/W)$ is 
the natural pro-algebraic group of the F-crystal $\H_\cris^{i}(X/W)$, it remains to apply Proposition \ref{s1}.

(iii) This 
is obtained by applying the derived functor $\otimes^L \Z/p^n$ 
to (\ref{n5})  and considering the associated long exact sequence,
using that $\Nyg{}\bH_\cris^{i}(X/W)$ and $\bH_\cris^{i}(X/W)$ are torsion-free.
\hfill $\Box$

\brem{\rm \label{str}
If $X$ is straight, then $\H^0(X,Z\Omega^1_X)\to \H^0(X,\Omega^1)$ is an isomorphism.
(Here by definition $Z\Omega^1_X=\Ker[d\colon \Omega^1_X\to\Omega^2_X]$.)
Then \cite[Thm.~II.5.14]{Ill} gives $\H^2(X,\Z_p(1))\cong\H^2_\cris(X/W)^{F=p}$, 
so we have $\H^2(X,\Z_p(1))[p^\infty]=0$. 
From Proposition \ref{fppfsplit} (ii) we deduce the surjectivity of the map
$$F/p-1:\Nyg{}\bH_\cris^1(X/W)\to\bH_\cris^1(X/W).$$ 
Thus the connected unipotent group attached to the F-crystal $\bH_\cris^1(X/W)$ is zero.
The exact sequence of Proposition \ref{fppfsplit} (iii) gives an isomorphism 
\begin{equation}
\bH^2(X,\mu_{p^n})\cong\Ker[\Nyg\bH_\cris^2(X/W)/p^n\xrightarrow{F/p-1}\bH_\cris^2(X/W)/p^n].
\label{mu}
\end{equation}
}
\erem

\bco \label{xx1}
Let $X$ be a straight variety over an algebraically closed field of characteristic $p>0$.
Then the quotient of $\Br(X)$ by its maximal divisible subgroup is the group of
$k$-points of the connected unipotent quasi-algebraic group associated to the F-crystal
$\H^2_\cris(X/W)$.
\eco

In this paper we are primarily interested in the case when $X=A$ is an abelian variety.

\bde \label{defU_A}
Let $A$ be an abelian variety over an algebraically closed field $k$ of characteristic $p>0$.
We denote by $\bU_A$ the connected unipotent quasi-algebraic group associated to the F-crystal
$\H^2_\cris(A/W)$.
\ede

In this case we have
$\H^2_\cris(A/W)\cong\wedge^2_W\H^1_\cris(A/W)$, where $\H^1_\cris(A/W)$ is the Dieudonn\'e
module of the $p$-divisible group $A[p^\infty]$.

The first statement of the following corollary is \cite[Cor.~5.2.4]{Gr}.

\bco \label{1.4}
Let $A$ be an abelian variety over an algebraically closed field $k$ of characteristic $p>0$, 
and let $i\geqslant 1$. The multiplication by $m$ map $[m]:A\to A$
acts on $\H^{i}(A,\Z_p(1))/\mathrm{tors}$ as multiplication by $m^i$, on 
$\H^i(X,\Z_p(1))[p^\infty]$ as multiplication by $m^{i-1}$, and on $\H^i_\fppf(A,\mu_{p^n})$,
for $i\leqslant 2$, and $\Br(A)$ as multiplication by $m^i$. 
\eco
{\em Proof.} The statements concerning $\H^{i}(A,\Z_p(1))$ follow from the exact sequence of 
Proposition \ref{fppfsplit} (ii) and the well-known fact that
$[m]$ acts as $m^i$ on $\H_\cris^{i}(A)\cong\wedge^i\H^1_\cris(A)$.

The statements about $\H^i_\fppf(A,\mu_{p^n})$ are well-known;
we give a proof for the convenience of the reader. 
The Kummer sequence gives $\H^0_\fppf(A,\mu_{p^n})=0$, so in degree 0 there is nothing to prove. It also gives a natural isomorphism
$\H^1_\fppf(A,\mu_{p^n})\cong\Pic(A)[p^n]$, which proves the statement in degree 1.
The statement in degree 2 follows from the functoriality of the injective map
$\H^2_\fppf(X,\mu_{p^n})\to\Nyg\H_\cris^2(X)/p^n$ recalled in Remark \ref{str}. See also
\cite[Cor.~5.7]{Yang}. The last statement follows from the Kummer exact sequence. \hfill $\Box$

\medskip

In general, the image of 
$\bH_\cris^2(A/W)\to\bH^3(A,\mu_{p^n})$ for an abelian variety $A$ can be non-trivial, for example when $A$ is the product of two supersingular elliptic curves \cite[Prop.~5.2.6]{Gr}. This shows that $[m]$ does not always act on $\H^3(A,\mu_{p^n})$ as multiplication by $m^3$.

\section{Unipotent groups from Dieudonn\'e modules}

\subsection{Unipotent groups attached to Dieudonn\'e modules}

Let $M$ and $N$ be Dieudonn\'e modules. Let us write $M\otimes N:=M\otimes_W N$.
The $W$-modules $M\otimes N$ and $\wedge^2 M$ are F-crystals
with the $\si$-linear action of $F$ given
by $F(x\otimes y)=Fx \otimes F y$ and $F(x\wedge y)=Fx\wedge Fy$.
A natural isomorphism of F-crystals 
\begin{equation}
\wedge^2 (M\oplus N)\cong (\wedge^2 M)\oplus (\wedge^2 N)\oplus (M\otimes N)
\label{mmm}
\end{equation}
gives a natural isomorphism of attached connected unipotent quasi-algebraic $k$-groups
(see Definition \ref{defU_L}):
\begin{equation}
\bU_{\wedge^2(M\oplus N)}\cong \bU_{\wedge^2 M}\times \bU_{\wedge^2 N}\times \bU_{M\otimes N}. \label{iso1}
\end{equation}

\brem {\rm
If $N\subset M$ is an inclusion of Dieudonné modules, we write $M\wedge N$ 
for the image of $M\otimes N\subset M\otimes M$ in $\wedge^2M$.}
\erem

\ble \label{5.1}
We have $\Nyg(M)=VM$, $\Nyg(M\otimes N)= VM\otimes N + M\otimes VN$,
and $\Nyg(\wedge^2 M)=M\wedge VM$. In particular, the canonical map
$\Nyg(M\otimes M)\to\Nyg(\wedge^2 M)$ is surjective.
\ele
{\em Proof.} Since $F$ is injective, we have $Fx=py$ if and only if $x=Vy$, giving
the first equality.

Let $m=\rk_W(M)$ and $n=\rk_W(N)$.
We can reduce modulo $p$ and consider $k$-vector spaces $M_1=M/p$ and 
$N_1=N/p$ of dimensions $m$ and $n$, respectively. Let
$m'=\dim_k(M_1[F])$ and $n'=\dim_k(N_1[F])$. Let $K$ be
the kernel of $F$ acting on $M_1\otimes_kN_1$. It is clear that
$$M_1[F]\otimes_k N_1+M_1\otimes_kN_1[F]\subset K.$$
This subspace has dimension $m'n+mn'-m'n'=mn-(m-m')(n-n')$. The dimension of
$F(M_1\otimes_k N_1)=FM_1\otimes_k FN_1$ is
$(m-m')(n-n')=mn-\dim(K)$, so the displayed inclusion is
an equality. Since $M_1[F]=VM_1$ and $N_1[F]=VN_1$,
we have $K=VM_1\otimes_kN_1+M_1\otimes_kVN_1$. 
The proof for $\wedge^2 M$ is similar. \hfill $\Box$

\medskip

Note that from $\Nyg(M)=VM$ we deduce $\bU_M=0$. Indeed, the map $1-V$ is surjective on $M$. 

\brem{\rm \label{e1}
For any Dieudonn\'e module $M$ we have
$$
\dim(\bU_{\wedge^2M}[p])\geqslant a(a-1)/2,$$
where $a=\dim_k(M/(FM+VM))$ 
is the $a$-number of $M$. To prove this we note that $F/p-1$
sends $\Nyg(\wedge^2 M)=M\wedge VM$ into $M\wedge (FM+VM)$. 
This gives a surjection $\bU_{\wedge^2M}\to\bG_{a,k}^{a(a-1)/2}$.
Comparing with (\ref{dim}) we see that this map is an isomorphism when
$a=g$, which happens exactly when $M$ is the Dieudonn\'e module of the superspecial abelian variety
$E^g$, where $E$ is a supersingular elliptic curve. 
}\erem

We note that $\Hom_W(M,N)$ has a natural structure of a pro-algebraic group, which is 
isomorphic to $\bW^{r_1r_2}$, where $r_1=\rk_W(M)$ and $r_2=\rk_W(N)$. This defines
a quasi-algebraic group structure on $\Hom_W(M/p^n,N/p^n)$
and hence also on $\Hom_\E(M/p^n,N/p^n)$ which is the subgroup of
$\Hom_W(M/p^n,N/p^n)$ invariant under the action of endomorphisms $F$ and $V$.

\bpr \label{prr}
Let $M$ and $N$ be Dieudonn\'e modules. For every $n\geqslant 1$
there is a natural isomorphism of quasi-algebraic groups 
$$\bU_{M^\vee\otimes N}[p^n]\cong\Hom_\E(M/p^n,N/p^n)/\left(\Hom_\E(M,N)/p^n\right).$$
\epr
{\em Proof.} Applying Corollary \ref{s3} with  $L=M^\vee\otimes N=\Hom_W(M,N)$ we obtain that 
$\bU_L[p^n]$ is naturally isomorphic to the cokernel of $\bL^{F=p}/p^n \to \bH_{L,n}$,
where by definition we have $\bH_{L,n}=\Ker[F/p-1\colon\Nyg(\bL)/p^n\to \bL/p^n]$.

The Frobenius $F$ acts on $L=\Hom_W(M,N)$ by sending
a homomorphism of $W$-modules $f\colon M\to N$ to $FfV$.
We have $f\in L^{F=p}$ if and only if $pf=FfV$,
which is equivalent to
$Vf=fV$, and this implies $Ff=fF$. We conclude that $L^{F=p}=\Hom_\E(M,N)$.

Next we note that $\Nyg(L)$ is the set of $f\in L$ such that $p$ divides $FfV$. 
Thus  $H_{L,n}=K_{L,n}/p^n\Nyg(L)$, where 
$$K_{L,n}:=\{f\in L \ \text{such that}  \ p | FfV(x)  \ \text{and}  \ p^n | FfV(x)/p-f(x)
\ \text{for any}  \ x\in M\}.$$
The conjunction of these two conditions is equivalent to 
$p^{n+1}|FfV(x)-pf(x)$, thus
$$K_{L,n}=\{f\in L \ \text{such that}  \ p^{n+1}\  \text{divides} \ FfV(x)-pf(x)
\ \text{for any}  \ x\in M\}.$$
Replacing $x$ by $Fx$ we deduce that
$Ff\equiv fF\bmod p^n$. Similarly, applying $V$ to both sides we get
$fV\equiv Vf\bmod p^n$. Thus reducing $f$ modulo $p^n$ gives a homomorphism 
$$\theta\colon K_{L,n}\to \Hom_\E(M/p^n,N/p^n).$$
By Corollary \ref{s3}, it remains to prove that $\theta$ is surjective with $\Ker(\theta)=p^n\Nyg(L)$.

It is clear that $p^n\Nyg(L)\subset \Ker(\theta)$. Conversely, let $f\in \Ker(\theta)$.
Then $f=p^ng$ for some $g\in\Hom_W(M,N)$ and we have that $p^{n+1}$ divides
$p^n(FgV-pg)$. Hence $p | FgV$ so that $g\in\Nyg(L)$. Thus $\Ker(\theta)=p^n\Nyg(L)$.

To prove that $\theta$ is surjective, take any $\bar f\in \Hom_\E(M/p^n,N/p^n)$ and lift it
to some $f\in \Hom_W(M,N)$. We need to find a $\phi\in\Hom_W(M,N)$ such that 
$f-p^n\phi\in K_{L,n}$, that is, such that $p^{n+1}|FfV-pf-p^n F\phi V$.

Since $Ff- fF\equiv 0\bmod p^n$, there is a $\si$-linear map
$g\colon M\to N$ such that $Ff- fF=p^ng$. This implies 
$FfV-pf=p^ngV$. Likewise, there is a $\si^{-1}$-linear map $h\colon M\to N$ such that
$fV-Vf=p^nh$. This implies $FfV-pf=p^n Fh$. The condition on $\phi$ can now be restated as
the condition that $p$ divides $gV-F\phi V=Fh-F\phi V$.

Since $pM\subset VM\subset M$, we can find $W$-bases $e_1,\ldots,e_n$ and $e'_1,\ldots,e'_n$
of $M$
such that $Ve_i=e'_i$ for $i=1,\ldots, m$ and $Ve_i=pe'_i$ for $i=m+1,\ldots,n$, for some 
$m\leqslant n$. We now define $\phi$ so that $\phi(e'_i)=h(e_i)$ for $i=1,\ldots, m$ and
$\phi(e'_i)=0$ for $i=m+1,\ldots,n$. Then $(Fh-F\phi V)(e_i)=0$ for $i=1,\ldots, m$.
For $i=m+1,\ldots,n$ we know that $p$ divides $Ve_i$, hence $p$ divides 
$(gV-F\phi V)(e_i)$. Thus the chosen $\phi$ satisfies the required condition. \hfill $\Box$

\bco \label{2.14}
Assume $p\neq 2$. Let $N$ be a Dieudonn\'e module. For every $n\geqslant 1$
there is a natural isomorphism of quasi-algebraic groups 
$$\bU_{\wedge^2 N}[p^n]\cong
\Hom_\E(N^\vee/p^n,N/p^n)^\skew/\left(\Hom_\E(N^\vee,N)^\skew/p^n\right),$$
where $\skew$ denotes the subgroup of anti-self-dual elements.
\eco
{\em Proof.} Let  $\iota$ be the involution permuting the factors of $N\otimes N$.
Since $p\neq 2$, we have isomorphisms
$$\wedge^2(N/p^n)\cong ((N/p^n)\otimes (N/p^n))^{\iota=-1}\cong
\Hom_W(N^\vee/p^n,N/p^n)^\skew$$
and compatible isomorphisms
$$\wedge^2N\cong (N\otimes N)^{\iota=-1}\cong \Hom_W(N^\vee,N)^\skew$$
via the reduction mod $p$ maps. Passing to the attached unipotent quasi-algebraic groups
we get an isomorphism
$\bU_{\wedge^2N}\cong \bU_{(N\otimes N)^{\iota=1}}\cong (\bU_{N\otimes N})^{\iota=1}$.
We thus obtain the desired isomorphism of unipotent quasi-algebraic
groups by taking the $\iota$-anti-invariant subgroups in the isomorphism of Proposition \ref{prr}
with $M=N$. \hfill $\Box$

\subsection{Behaviour of $\bU_{\wedge^2M}$ and $\bU_{M\otimes N}$ under isogeny}

Let $K=W[1/p]$ be the field of fractions of $W$. 
We call Dieudonn\'e modules $M$ and $M'$ {\em isogenous} if the $\E$-modules
$M\otimes_WK$ and $M'\otimes_WK$ are isomorphic.

\bpr \label{iso-inv}
The integers $\dim(\bU_{\wedge^2M})$ and $\dim(\bU_{M\otimes N})$ depend only on the isogeny classes of $M$ and $N$.
\epr
{\em Proof.} By (\ref{iso1}) it is enough to prove the statement for $\bU_{\wedge^2M}$. 

Suppose that there is an exact sequence of $\E$-modules
\begin{equation}
0\to M\to M'\to C\to 0, \label{ext}
\end{equation}
where $p^n C=0$ for some $n$. In particular, we have $\rk_WM=\rk_W M'$.
We can consider (\ref{ext}) as an exact sequence in the abelian category ${\mathcal P}$ of
pro-algebraic $k$-groups, with the maps respecting the action of $F$ and $V$ by $k$-endomorphisms.
Applying the snake lemma to the commutative diagram
obtained from the action of $V$ on $M$, $M'$, $K$,
we get an exact sequence of quasi-algebraic groups
$$0\to C[V]\to M/VM\to M'/VM'\to C/VC\to 0.$$
Counting dimensions we see that $\dim_k(M/VM)=\dim_k(M'/VM')$.
This implies 
$$\dim_k( \wedge^2 M/(M\wedge VM))=\dim_k( \wedge^2 M'/(M'\wedge VM')).$$
Defining the quasi-algebraic groups
$$C_1=(M'\wedge VM')/(M\wedge VM), \quad C_2=\wedge^2 M'/\wedge^2 M,$$
we deduce that $\dim(C_1)=\dim(C_2)$. We have the following commutative diagram
of pro-algebraic groups
$$\xymatrix{0\ar[r]&M\wedge VM\ar[r]\ar[d]_{F/p-1}&M'\wedge VM'\ar[r]\ar[d]_{F/p-1}&C_1\ar[r]\ar[d]_\phi&0\\
0\ar[r]&\wedge^2 M\ar[r]&\wedge^2 M'\ar[r]&C_2\ar[r]&0}$$
where the rows are exact and the map $\phi$ is induced by two other vertical maps.
The snake lemma gives an exact sequence
$$0\to (\wedge^2M)^{F=p}\to  (\wedge^2M')^{F=p}\to \Ker(\phi)\to \bU_{\wedge^2M}\to \bU_{\wedge^2M'}\to \Coker(\phi)\to 0.$$
Now $\dim(C_1)=\dim(C_2)$ implies that $\dim(\Ker(\phi))=\dim(\Coker(\phi))$. The quotient of $(\wedge^2M')^{F=p}$ by $(\wedge^2M)^{F=p}$ is a finite group. Thus
counting dimensions of the terms of the exact sequence 
we obtain $\dim(\bU_{\wedge^2M})=\dim(\bU_{\wedge^2M'})$. \hfill $\Box$

\medskip

Unlike the dimension, the $p$-exponent of $\bU_{\wedge^2M}$ changes under isogeny.
In certain cases, one can control this change.

\bpr \label{livia}
Let $M$ and $M'$ be Dieudonn\'e modules that fit into an exact sequence of $\E$-modules 
$$0\to M\to M'\to C\to 0.$$

{\rm (i)} If $VC=0$, then the kernel of the natural map 
$\bU_{\wedge^2M}(k)\to \bU_{\wedge^2M'}(k)$ is annihilated by $p$.

{\rm (ii)} If $FC=0$, then the cokernel of 
$\bU_{\wedge^2M}(k)\to \bU_{\wedge^2M'}(k)$ is annihilated by $p$. 

{\rm (iii)} If $VC=FC=0$, then we have an exact sequence
$$1\to\bG_{a,k}^n\to \bU_{\wedge^2M}\to \bU_{\wedge^2M'}\to \bG_{a,k}^n\to 1$$
for some $n\geqslant d(d-1)/2$, where $d=\dim(C)$. In particular,
the $p$-exponents of $\bU_{\wedge^2M}$ and $\bU_{\wedge^2M'}$ differ at most by~$1$.
\epr
{\em Proof.} 
We can think of $\wedge^2 M$, $M\wedge VM$, $M'\wedge VM'$ as subcrystals of $\wedge^2M'$.

In the notation of the previous proof we show that $p\Ker(\phi)=0$ if $VC=0$,
and $p\Coker(\phi)=0$ if $FC=0$. This will prove (i) and (ii).

To prove the first statement, let $a\in M'\wedge VM'$ be such that 
$(F/p)(a)-a\in \wedge^2 M$. We need to prove that $pa\in M\wedge VM$.  
The condition $VC=0$ implies that $VM'\subset M$, and this implies
$pM'\wedge VM'=VFM'\wedge VM'\subset \wedge^2 M$. In particular, we have 
$pa\in \wedge^2 M$, and hence $F(pa)\in p\wedge^2 M$. This says that 
$pa\in\Nyg(\wedge^2 M)$. By Lemma \ref{5.1} we have $\Nyg(\wedge^2 M)=M\wedge VM$,
so $pa\in M\wedge VM$, as desired. 

To prove the second statement, let $b\in\wedge^2M'$. 
It is clear that $pb\in\Nyg(\wedge^2M')=M'\wedge VM'$. We can write
$pb=Fb-(F/p-1)(pb)$. The condition $FC=0$ implies that $FM'\subset M$, and this implies
$F(\wedge^2M')\subset \wedge^2 FM'\subset\wedge^2M$. Hence the image of $pb$ in 
$\Coker(\phi)$ is zero, as required. 

(iii) If $VC=FC=0$, then $F/p-1$ sends $M'\wedge VM'$ into $M\wedge M'\subset\wedge^2M'$.
This gives
$n\geqslant\dim_k\left(\wedge^2M'/(M\wedge M')\right)=d(d-1)/2$. 
\hfill $\Box$

\subsection{The dimensions of $\bU_{\wedge^2 M}$ and $\bU_{M\otimes N}$}

We now calculate $\dim(\bU_{\wedge^2 M})$ in terms of the Newton polygon of $M$,
using only the classical Dieudonn\'e theory and the invariance of 
$\dim(\bU_{\wedge^2 M})$ under the isogeny of $M$. 

\bthe \label{th-dim}
Let $M$ be a Dieudonn\'e module with slopes $\lambda_i=c_i/(c_i+d_i)$
of multiplicity $c_i+d_i$, where $c_i$ and $d_i$
are coprime non-negative integers, for $i=1,\ldots,n$ (allowing repetitions among the pairs
$(c_i,d_i)$). 
Let $m_i\coloneqq\min(c_i,d_i)$. Then we have 
$$\dim(\bU_{\wedge^2 M})=\sum_{i=1}^n m_i(m_i-1)/2 +\sum_{1\leq i<j\leq n}\min(c_ic_j,d_id_j).$$
\ethe
{\em Proof.}
By Proposition \ref{iso-inv}, $\dim(\bU_{\wedge^2 M})$
is an isogeny invariant of $M$, so we can replace $M$ by a direct sum
of arbitrary isoclinic Dieudonn\'e modules $M_i$ of rank $c_i+d_i$ and slope $c_i/(c_i+d_i)$.
Using (\ref{mmm}) we obtain
\begin{equation}
\dim(\bU_{\wedge^2 M})=\sum_{i=1}^n\dim(\bU_{\wedge^2 M_i})+
\sum_{1\leq i<j\leq n}\dim(\bU_{M_i\otimes M_j}). \label{nnn}
\end{equation}

Theorem \ref{th-dim} is a consequence of Proposition \ref{5.8} and Corollary \ref{5.10} below.

\bpr \label{5.8}
Let $M_i$ for $i=1,2$ be isoclinic Dieudonn\'e modules of rank $c_i+d_i$ and slope 
$c_i/(c_i+d_i)$, where $c_i,d_i$ are coprime and non-negative. Then we have
$$\dim(\bU_{M_1\otimes M_2})=\min(c_1c_2,d_1d_2).$$ 
\epr
{\em Proof.} For integers $m\leqslant n$ we denote by $[m,n]$ the set of integers $m,m+1,\ldots, n$.

By Proposition \ref{iso-inv} we can take $M_1=M_{c_1,d_1}$. 
This is a free $W$-module with basis $e_i$ for $i\in [1,c_1+d_1]$ such that
$Fe_i=e_{i+c_1}$ and $Ve_i=e_{i+d_1}$ with the convention that $e_{i+c_1+d_1}=pe_i$.
Let us take $M_2=M_{c_2,d_2}$ and denote the basis vectors of $M_{c_2,d_2}$ by
$f_j$ for $j\in [1,c_2+d_2]$. 
Write $H:=M_{c_1,d_1}\otimes_W M_{c_2,d_2}$. This is a free $W$-module with basis 
$v_{i,j}=e_i\otimes f_j$. Here $$(i,j)\in S_0 \cup S_1\cup S_2\cup S_3,$$
where $S_0=[1,d_1]\times[1,d_2]$, $S_1=[d_1+1,c_1+d_1]\times[1,d_2]$, 
$S_2=[1,d_1]\times[d_2+1,c_2+d_2]$, $S_3=[d_1+1,c_1+d_1]\times[d_2+1,c_2+d_2]$.
We note that $Fv_{i,j}=v_{i+c_1,j+c_2}$ with the convention that 
$v_{i+c_1+d_1,j}=v_{i,j+c_2+d_2}=pv_{i,j}$.

To prove the proposition, it is enough to show that $\bU_H$ is the product of $\min(c_1c_2,d_1d_2)$ copies of $\bG_{a,k}$. The proof of this consists of a series of claims.

\medskip

\noindent {\bf Claim 1}. $\Nyg(H)$ is the $W$-submodule of $H$ with basis
$pv_{i,j}$ for $(i,j)\in S_0$ and $v_{i,j}$ for $(i,j)\in S_1\cup S_2\cup S_3$. This is clear.

\medskip

\noindent {\bf Claim 2}. Let $\tau$ be an orbit of $F$ acting on $S$ by translation by $(c_1,c_2)$ subject to the convention mentioned above, and let $M_\tau$ be the $W$-submodule of 
$H$ generated by $v_{i,j}$ for $(i,j)\in\tau$. The cokernel of
\begin{equation}
F/p-1\colon\Nyg(M_{c_1,d_1}\otimes M_{c_2,d_2})\to M_{c_1,d_1}\otimes M_{c_2,d_2}
\label{Nyg-map}
\end{equation}
is the direct sum of 
$$\bU_\tau:=\Coker[F/p-1\colon\Nyg(M_\tau)\to M_\tau],$$
where $\Nyg(M_\tau)=\{x\in M_\tau| p \  \text{divides} \  Fx\}$, over all orbits $\tau$. This is also clear.

\smallskip

We now fix an orbit $\tau$ of $F$ in $S$. We can think of $\tau$ as a cyclic graph with an orientation
and refer to the elements of $\tau$ as vertices.

\medskip

\noindent {\bf Claim 3}. Assume that $d_1d_2\leqslant c_1c_2$.
If $x,y\in\tau\cap S_0$, then there exists a vertex $z\in \tau\cap S_3$ between $x$ and $y$
in the direction from $x$ to $y$. Indeed,
in the opposite case we have $z\in S_1\cup S_2$ for all $z\in \tau$ between $x$ and $y$.
Let us show that this leads to a contradiction. Write $x=(x_1,x_2)$ and $y=(y_1,y_2)$.
Suppose that there are exactly $t_1$ vertices of $S_1$ and $t_2$ vertices of $S_2$
in $\tau$ between $x$ and $y$ in the direction from $x$ to $y$. Using the 
description of the action of $F$ we obtain
\begin{align*}
y_1-x_1&=(t_1+t_2+1)c_1-t_1(c_1+d_1),\\
y_2-x_2&=(t_1+t_2+1)c_2-t_2(c_2+d_2).
\end{align*}
Since $|y_1-x_1|<d_1$ and $|y_2-x_2|<d_2$, we get
$$c_1(t_2+1)<d_1(t_1+1), \quad c_2(t_1+1)<d_2(t_2+1).$$
Multiplying these inequalities we obtain $c_1c_2<d_1d_2$, which is a contradiction.

A similar proof shows that if $d_1d_2\geqslant c_1c_2$, then $\tau$ contains 
a vertex in $S_0$ between any two vertices in $S_3$.

\medskip

\noindent {\bf Claim 4}. Let us use the term {\em chain} 
for a connected and simply connected subgraph of $\tau$.
A chain will be called {\em relevant} if all its vertices are in $S_1\cup S_2$ except the last one,
which is in $S_0$; moreover, the vertex immediately preceeding the first vertex of the chain
must be in $S_3$. 
The elements of relevant chains in $\tau$ will be called {\em relevant}.
Let $L_\tau$ be the $W$-submodule of $M_\tau$ generated by $pv_x$ if $x\in\tau$ is relevant 
and by $v_x$ if $x\in \tau$ is not relevant. We claim that $L_\tau$ is contained 
in $(F/p-1)(\Nyg(M_\tau))$.

To prove this we define a linear map $\psi\colon L_\tau\to L_\tau$ by the rule: if $x$ is relevant,
then $\psi(pv_x)=pv_x$; if $x\in S_0$ is not relevant, then $\psi(v_x)=pv_x$; and if $x\notin S_0$
is not relevant, then $\psi(v_x)=v_x$. Next, define a $\si$-linear map $\phi\colon L_\tau\to L_\tau$ 
by the rule: if $x$ is relevant, then $\phi(pv_x)=F(v_x)$; if  $x\in S_0$ is not relevant, 
then $\phi(v_x)=F(v_x)$; and if $x\notin S_0$ is not relevant, then $\phi(v_x)=F(v_x)/p$.
It is easy to check that $\phi(L_\tau)=L_\tau$, and by semilinear algebra
this implies that $(\phi-\psi)(L_\tau)=L_\tau$. A direct verification shows that
$(\phi-\psi)(L_\tau)$ is contained in $(F/p-1)(\Nyg(M_\tau))$, so the claim is proved.

\medskip

\noindent {\bf Claim 5}: $\bU_\tau\simeq\bG_{a,k}^n$, where $n$ is the number of relevant chains in $\tau$.

Let us prove this. The $k$-vector space $M_\tau/L_\tau$ has a basis consisting of the images of
$v_x$ for all relevant $x\in\tau$. We can think of $M_\tau/L_\tau$ as a quasi-algebraic group
given by a product of copies of $\bG_{a,k}$ numbered by the relevant vertices in $\tau$.

Suppose that $v_1,\ldots,v_r$ is a relevant chain in $\tau$.
For any $i=1,\ldots,r-1$, we have $v_i\in \Nyg(M_\tau)$,
thus $\sum_{i=1}^{r-1}(F/p-1)(W v_i)$
is contained in $(F/p-1)(\Nyg(M_\tau))$. 
This gives an injective map of quasi-algebraic groups
$\bG_{a,k}^{r-1}\to\bG_{a,k}^r$, so its cokernel is isomorphic to $\bG_{a,k}$. In fact, $(F/p-1)(\Nyg(M_\tau))$ is the sum of $L_\tau$ and the subgroups
$\sum_{i=1}^{r-1}(F/p-1)(W v_i)$ for all relevant chains in $\tau$.
Therefore, each relevant chain in $\tau$ contributes a copy of $\bG_{a,k}$ to $\bU_\tau$. This proves the claim.

\medskip

\noindent{\em End of proof of Proposition} \ref{5.8}.
From Claim 2 we easily deduce that the number of relevant chains in
all orbits is $\min(c_1c_2,d_1d_2)$. Indeed, if $d_1d_2\leqslant c_1c_2$, then associating to a relevant
chain to its last vertex defines a bijection between the set of relevant chains and $S_0$.
Similarly, if $d_1d_2\geqslant c_1c_2$, then associating to a relevant chain the vertex 
immediately preceeding its first vertex defines a bijection between the set of relevant chains 
and $S_3$. Applying Claim 5 we prove the proposition. \hfill $\Box$

\brem{\rm It is likely that
this proposition can also be deduced from \cite[Examples 2.3.1, 2.3.2]{Vas08}.
}\erem

\bco \label{5.10}
Let $M$ be an isoclinic Dieudonn\'e module of rank $c+d$ and slope 
$c/(c+d)$, where $c,d$ are coprime and non-negative. 
Then we have 
$$\dim(\bU_{\wedge^2 M})=m(m-1)/2,$$ 
where $m=\min(c,d)$.
\eco
{\em Proof.} The $W$-module $\wedge^2 M$ is the quotient of $M\otimes M$ by 
the $W$-submodule $K$ generated by $x\otimes x$, for $x\in M$, that is, 
in the notation of the previous proof, by $v_{i,i}$ for all $i$ and by $v_{i,j}+v_{j,i}$ for all $i\neq j$.
By Lemma \ref{5.1} the natural map
$\Nyg(M\otimes M)\to\Nyg(\wedge^2 M)$ is surjective, so we have a commutative
diagram
$$\xymatrix{&&\Nyg(M\otimes M)\ar[r]\ar[d]^{F/p-1}&
\Nyg(\wedge^2 M)\ar[r]\ar[d]^{F/p-1}&0\\
0\ar[r]&K\ar[r]&M\otimes M\ar[r]&\wedge^2 M\ar[r]&0}$$
Thus $\bU_{\wedge^2 M}(k)$ is the quotient of $\bU_{M\otimes M}(k)$ by the image of $K$.

We now assume $d\leqslant c$, the proof in the opposite case being entirely similar.
By the proof of Proposition \ref{5.8}, for $M=M_{c,d}$, the group $\bU_{M\otimes M}$
is the product of copies of $\bG_{a,k}$ numbered by the elements of $S_0=[1,d]\times[1,d]$.
(Indeed, each relevant chain contributes a copy of $\bG_{a,k}$, but for every $x\in S_0$ there is exactly
one relevant chain that ends in $x$.) The involution swapping the factors
of $M\otimes M$ acts on the complement to the diagonal in $[1,d]\times[1,d]$ without fixed points,
thus the image of $K$ in $\bU_{M\otimes M}(k)\cong k^{d^2}$ is a vector subspace of 
dimension $d+d(d-1)/2$. Hence  $\bU_{\wedge^2 M}\cong \bG_{a,k}^{d(d-1)/2}$. By Proposition \ref{iso-inv} this proves the corollary.
\hfill $\Box$

\medskip

This finishes the proof of Theorem \ref{th-dim}.

\bco \label{dim 0}
Let $M$ be a Dieudonn\'e module.  Then $\bU_{\wedge^2 M}=0$ if and only if
the slopes of $M$ are $0$ (with any non-negative multiplicity), 
$1$ (with any non-negative multiplicity), and at most one slope $1/(1+n)$ or $n/(1+n)$
with multiplicity $1+n$.

\eco
{\em Proof.} The formula of Theorem \ref{th-dim} gives $\dim(\bU_{\wedge^2 M})=0$ if and only if
$m_i\in\{0,1\}$ for all $i$, and $m_i=1$ for at most one subscript $i$. 
This implies the statement. \hfill $\Box$

\bco \label{d3}
Let $M$ and $N$ be a Dieudonn\'e modules. Then $\bU_{M\otimes N}=0$ if and only if
at least one of $M$ and $N$ has all of its slopes $0$ (with any non-negative multiplicity) and
$1$ (with any non-negative multiplicity).
\eco
{\em Proof.} This immediately follows from Proposition \ref{5.8}. \hfill $\Box$

\section{Abelian varieties} 

\subsection{The dimension of $\bU_A$} \label{S3.1}

Using the results of the previous section it is now easy to deduce a formula for $\dim(\bU_A)$,
where $A$ is an abelian variety over $k$.

\bthe \label{3.1}
Let $A$ be an abelian variety of dimension $g$ over $k$. Then we have
\begin{equation}
\dim(\bU_A)=\frac{g(g-1)}{2}-\sum_{0\leqslant \lambda<1}(1-\lambda)m_\lambda, \label{u1}
\end{equation}
where $m_\lambda$ is the multiplicity of slope $\lambda$ in $\H^2_\cris(A/W)$.
\ethe
{\em Proof.} The contravariant Dieudonn\'e module of $A$ is $M:=\H^1_\cris(A/W)$,
and there is a natural isomorphism $\H^2_\cris(A/W)\cong\wedge^2 M$.
Thus we have $\bU_A=\bU_{\wedge^2M}$. By Proposition \ref{iso-inv} we can assume
without loss of generality that $M$ is the direct sum of
isoclinic Dieudonn\'e modules $M_{c,d}$ of rank $c+d$ and slope $c/(c+d)$,
where $c,d$ are non-negative coprime integers.
Since $M$ comes from an abelian variety, this decomposition is invariant under the symmetry
swapping $c$ and $d$. Thus we can write
$$M=M_0\oplus M_1\oplus\ldots\oplus M_n,$$
where $M_0\cong M_{1,1}^{\oplus m}$ for some $m\geqslant 0$, and 
$M_i=M_{c_i,d_i}\oplus M_{d_i,c_i}$ for $i=1,\ldots,n$,
where $c_i, d_i$ are non-negative and coprime with $c_i<d_i$.
We have
\begin{equation}
\dim(\bU_{\wedge^2 M})=\sum_{i=0}^n \dim(\bU_{\wedge^2 M_i})+
\sum_{0\leqslant i<j\leqslant n} \dim(\bU_{M_i\otimes M_j}). \label{u2}
\end{equation}
We can compute this using Proposition \ref{5.8} and Corollary \ref{5.10}.
Let us show by induction on $n$ that the answer agrees with formula (\ref{u1}).

This is true for $n=0$. 
Indeed, in this case $g=m$ and $\wedge^2 M=\wedge^2M_0$ is isoclinic of slope 1, so
(\ref{u1}) gives $g(g-1)/2$. Applying Proposition \ref{5.8} and Corollary \ref{5.10}
we check that (\ref{u2}) gives the same answer.

For $i\geqslant 1$, an application of  Proposition \ref{5.8} and Corollary \ref{5.10}
gives $\dim(\bU_{\wedge^2 M_i})=c_i(c_i-1)+c_id_i$. This agrees with 
formula (\ref{u1}) that gives
$$\frac{(c_i+d_i)(c_i+d_i-1)}{2}-\left(1-\frac{2c_i}{c_i+d_i}\right)\frac{(c_i+d_i)(c_i+d_i-1)}{2}=
c_i(c_i+d_i-1).$$

Let us assume now that the two formulae agree for $M$ 
and show that they also agree for $M\oplus M_{n+1}$. 
Let $g=m+\sum_{i=1}^n(c_i+d_i)$ be one half of the rank of $M$.
Write $M_{n+1}=M_{c,d}\oplus M_{d,c}$, where $c<d$. 
By the above, it remains to check that 
$$\dim(\bU_{M\otimes M_{n+1}})=g(c+d)-\sum_{0\leqslant\lambda<1}(1-\lambda)m_\lambda,$$
where $\lambda$ is a slope of $M\otimes M_{n+1}$. It is enough to check this
for $M=M_{1,1}$ and for $M=M_i$, where $i=1,\ldots,n$.

Let $N=M_{1,1}\otimes (M_{c,d}\oplus M_{d,c})$. We need to check that
$$\dim(\bU_N)=(c+d)-\sum_{0\leqslant\lambda<1}(1-\lambda)m_\lambda,$$
where $\lambda$ is a slope of $N$. Calculating slopes and their multiplicities, 
we obtain that the right hand side equals $2c$, and the same is obtained by applying 
Proposition \ref{5.8} and Corollary \ref{5.10}.

Now let $N=(M_{a,b}\oplus M_{b,a})\otimes (M_{c,d}\oplus M_{d,c})$, where
$a,b$ are coprime non-negative integers such that $a<b$. We need to check that
$$\dim(\bU_N)=(a+b)(c+d)-\sum_{0\leqslant\lambda<1}(1-\lambda)m_\lambda,$$
where $\lambda$ is a slope of $N$. Without loss of generality we can assume that
$a/b\geqslant c/d$. Calculating slopes and their multiplicities, we obtain that the right hand side of
the last displayed formula gives $2c(a+b)$, and the same is obtained by applying 
Proposition \ref{5.8} and Corollary \ref{5.10}. This finishes the proof. \hfill $\Box$

\brem{\rm \label{ferry}
(i) Formula (\ref{u1}) implies that $\dim(\bU_A)\leqslant g(g-1)/2$ with equality if and only if $g=1$ or
$A$ is supersingular. Indeed, the slopes of $\H^2_\cris(A/W)$ are invariant under the 
reflection $\lambda\mapsto 2-\lambda$, so the equality is attained if and only if
$\H^2_\cris(A/W)$ is isoclinic of slope $1$. For $g\geqslant 2$
this implies that $\H^1_\cris(A/W)$ is isoclinic of slope $1/2$.

(ii) We have $\bU_A=0$ if and only if $A$ is ordinary (i.e.,
the slopes of $\H^1_\cris(A/W)$ are $0$ and $1$, each with multiplicity $g$) or almost ordinary
(i.e., the slopes of $\H^1_\cris(A/W)$ are $0$ and $1$, each with multiplicity $g-1$, and $1/2$
with multiplicity $2$). This follows from Corollary \ref{dim 0} or from the proof of Theorem \ref{3.1}.

(iii) We have $\bU_A=\bG_{a,k}$ if and only if the slopes of $\H^1_\cris(A/W)$ are $0$ and $1$, 
each with multiplicity $g-2$, and $1/2$ with multiplicity 4. This can be deduced from the
dimension formula of Theorem \ref{th-dim}.
}
\erem

\bco \label{cor3.3}
The formal Brauer group of $A$ is smooth.
\eco
{\em Proof.} By the dimension formula (\ref{u1}), this follows from Proposition \ref{smbrhat} below. \hfill$\Box$

\subsection{Relation to the formal Brauer group} \label{formal}

Let $X$ be a smooth proper variety over $k$, and let $\bU_X$ be the connected component
of the quasi-algebraic group over $k$ whose group of $k$-points is $\H^3(X,\Z_p(1))[p^\infty]$ 
(as discussed in the introduction). 

In this section we relate $\bU_X$ to the formal Brauer group of $X$.
The formal Brauer group $\Brhat(X)$ is the functor parameterising deformations of the trivial class in 
$\Br(X)$. For an Artinian $k$-algebra $R$ one defines
\begin{equation}
    \Brhat(X):=\Ker[\Br(X_R)\to\Br(X_{R_{\textnormal{red}}})],
\end{equation}
where $X_R=X\times_k\Spec(R)$. A detailed exposition of the formal Brauer group and its properties can be found in \cite{artinmazur}.

\bpr\label{brhat}
Let $X$ be a smooth proper variety over $k$. The following properties hold.

$(1)$ If the Picard scheme ${\bf Pic}_{X/k}$ is smooth, then $\Brhat(X)$ is representable
by a connected formal group of finite type with tangent space $\H^2(X,\mathcal{O}_X)$.

$(2)$ If $\widehat{\Br}(X)$ is representable, then its Cartier module of typical curves 
is isomorphic to $\H^2(X,W\mathcal{O}_X)$.

$(3)$ If $\widehat{\Br}(X)$ is representable and $\H_\cris^{3}(X/W)$ is torsion-free, then 
$\Brhat(X)$ is smooth.

\epr
{\em Proof.} Parts (1) and (2) are proved in \cite{artinmazur}, see 
Proposition II.1.8, Corollaries II.2.4 and II.4.3. Part (3) is a consequence of work of Ekedahl \cite{Eke}. For a different proof, see \cite[Cor.~3.17]{Gr2}. \hfill$\Box$

\medskip

Suppose that ${\bf Pic}_{X/k}$ is smooth, so that $\Brhat(X)$ is representable,
 and let $\Brhat(X)_\red$ be the reduced subgroup of $\Brhat(X)$.
Recall \cite[(28.3.2)]{Ha}
that a smooth connected formal group of finite type is an extension of a smooth connected formal group of finite height by a smooth unipotent formal group, in a unique way. 
Let us denote by $\widehat{U}_X$ the unipotent part of $\Brhat(X)_\red$, and by $\widehat{P}_X$ the maximal finite height quotient of $\Brhat(X)_\red$. 
By construction, $\widehat{U}_X$ is the reduced formal subgroup of $\Brhat(X)[p^n]$ for large enough $n$.

Let $\pi\colon X\to\Spec(k)$ be the structure morphism.
By a theorem of Bragg and Olsson \cite[Cor.~1.6]{braggolsson} we know that
$R^2\pi_*\mu_{p^n,X}$ is represented by a group $k$-scheme. Let us denote
by ${\mathcal U}_{X,n}$ its reduced connected component. The exact sequence
$$1\to \mu_{p^n,X}\to\mu_{p^{n+1},X}\to \mu_{p,X}\to 1$$
and the natural isomorphism $R^1\pi_*\mu_{p^n,X}\cong {\bf Pic}_{X/k}[p^n]$
imply that the induced map $R^2\pi_*\mu_{p^n,X}\to R^2\pi_*\mu_{p^{n+1},X}$
is injective. Define  ${\mathcal U}_X:=\varinjlim {\mathcal U}_{X,n}$. The dimension of
${\mathcal U}_{X,n}$ stabilises, so ${\mathcal U}_X\cong{\mathcal U}_{X,n}$ for
$n$ large enough.

\ble\label{compl}
The formal group $\widehat{U}_X$ is canonically isomorphic to the formal 
completion of ${\mathcal U}_X$ at the identity.
\ele
{\em Proof.} A more general statement is proved in \cite[version~1, Prop.~10.7]{braggolsson}. 
Recall that the formal completion functor maps a presheaf 
$\mathcal{F}$ over $\Spec(k)_\fppf$ to the functor
\begin{equation*}
    \widehat{\mathcal{F}}:R\mapsto\Ker[\mathcal{F}(R)\to\mathcal{F}(R_{\textnormal{red}})]
\end{equation*}
from the category of Artinian $k$-algebras to the category of abelian groups. Since the projection $R\to R_{\textnormal{red}}$ has a canonical section, the formal completion functor is exact. Pushing forward the Kummer sequence along $\pi:X\to\Spec(k)$ we obtain the short exact sequence of fppf sheaves
$$   0\to R^1\pi_*\mathbb{G}_{m,X}/p^n\to R^2\pi_*\mu_{p^n,X}\to R^2\pi_*\mathbb{G}_{m,X}[p^n]\to 0,$$
where $ R^1\pi_*\mathbb{G}_{m,X}/p^n$ is representable by an étale group scheme, by the theory of the Picard scheme \cite[p.~121]{CTS21}. Therefore, taking the formal completion we get an isomorphism
\begin{equation}\label{isocompl}
    (R^2\pi_*\mu_{p^n,X})^{\wedge}\cong (R^2\pi_*\mathbb{G}_{m,X}[p^n])^{\wedge}.
\end{equation}
Notice that the right-hand group is by definition $\Brhat(X)[p^n]$. 
Take $n$ large enough so that $\widehat{U}_X$ is the reduced formal subgroup of $\Brhat(X)[p^n]$, and such that ${\mathcal U}_X={\mathcal U}_{X,n}$. 
Then \eqref{isocompl} gives a desired isomorphism $({\mathcal U}_X)^\wedge\cong\widehat{U}_X$. \hfill$\Box$

\medskip

By the Kummer sequence, the quasi-algebraic group $({\mathcal U}_X)^\perf$ attached to 
${\mathcal U}_X$ is isogenous to $\bU_X$. From Lemma \ref{compl} we obtain
\begin{equation} \label{dimdim}
 \dim(\bU_X)=\dim({\mathcal U}_X)=\dim(\widehat{U}_X)=
\dim(\Brhat(X)_\red)-\dim(\widehat{P}_X).
\end{equation}
Let us compute $\dim(\widehat{P}_X)$. By Proposition \ref{brhat} (2), the
Cartier module of typical curves of $\Brhat(X)_\red$
 is isomorphic to $\H^2(X,W\mathcal{O}_X)$.
It follows that, after inverting $p$, the covariant Dieudonné module of $\widehat{P}_X$
 is isomorphic to $\H^2(X,W\mathcal{O}_X)[1/p]$, hence it is isogenous to 
$\H^2_\cris(X/W)_{[0,1)}$ by the theory of the slope spectral sequence. It is then well-known
(see, e.g. \cite[Proof of Cor.~III.3.4]{artinmazur}) that
\begin{equation} \label{dimP}
    \dim(\widehat{P}_X)=\sum_{0\leqslant\lambda<1}m_{\lambda}(1-\lambda),
\end{equation}
where $\lambda$ is a slope of $\H^2_\cris(X/W)$ of multiplicity $m_\lambda$, see Lemma \ref{dimgrp} below. 

\ble\label{dimgrp}
Let $G$ be a $p$-divisible group. Let $M(G)$ be the covariant Dieudonné module of $G$.
We have
\begin{equation*}
    \dim(G)=\sum_{\lambda}m_{\lambda}(1-\lambda),
\end{equation*}
where $\lambda$ is a slope of $M(G)$ of multiplicity $m_{\lambda}$.
\ele
{\em Proof.} The dimension of a $p$-divisible group is an isogeny invariant, so we may assume that $M(G)$ is the simple isoclinic Dieudonn\'e module $M_{c,d}$ as in the proof of Proposition \ref{5.8}.
Thus $\lambda=c/(c+d)$ and $m_\lambda=c+d$. We can easily compute
$\dim(G)=\dim_k(M(G)/VM(G))=d=m_\lambda(1-\lambda)$. 
\hfill$\Box$


\bpr\label{smbrhat}
Let $X$ be a proper smooth variety over $k$ such that the Picard scheme
${\bf Pic}_{X/k}$ is smooth.  Then the formal Brauer group $\Brhat(X)$ is smooth if and only if 
\begin{equation*}
    h^{0,2}(X)=\dim(\bU_X)+\sum_{\lambda}m_{\lambda}(1-\lambda),
\end{equation*}
where $m_{\lambda}$ is the multiplicity of the slope $\lambda$ in $\H^2_\cris(X/W)_{[0,1)}$. 
\epr
{\em Proof.} By Proposition \ref{brhat} (1) the tangent space of $\Brhat(X)$ 
has dimension $h^{0,2}(X)=\dim(\H^2(X,\mathcal{O}_X))$. Thus
$\Brhat(X)$ is smooth if and only if $\dim(\Brhat(X)_\red)=h^{0,2}(X)$. 
It remains to apply (\ref{dimdim}) and (\ref{dimP}). \hfill $\Box$

\medskip

This gives another proof of Theorem \ref{3.1} assuming smoothness of the formal Brauer group
of abelian varieties.

When $X$ is a surface such that ${\bf Pic}_{X/k}$ is smooth, then $\Brhat(X)$ is smooth \cite[p.~190]{artinmazur}.
In this case Proposition \ref{smbrhat} follows from \cite[Remark 7.4]{Milne75}.

\subsection{The $p$-exponent of $\bU_A$ is at most $g-1$} 

In this section we prove the following theorem. The non-supersingular case uses Corollary \ref{2.14}
and thus requires the assumption $p\neq 2$.

\bthe \label{ss}
Let $A$ be an abelian variety of dimension $g\geqslant 1$ over an algebraically closed field $k$
of positive characteristic $p\neq 2$. Then $p^{g-1}\bU_A=0$.
\ethe

In Section \ref{exa} we show that the bound of Theorem \ref{ss} is best possible.

\bco \label{f.g.}
Let $k_0$ be a field finitely generated over $\F_p$, where $p\neq 2$. 
Let $k$ be an algebraic closure of $k_0$.
Let $A$ be an abelian variety over $k_0$ of dimension $g\geqslant 1$.
Then the transcendental Brauer group $\Im[\Br(A)\to\Br(A_k)]$ is a direct sum of 
a finite group and an abelian group annihilated by $p^{g-1}$. 
\eco
{\em Proof.} This follows from Theorem \ref{ss} and a theorem of D'Addezio
\cite[Thm.~5.2]{D'A}, see also \cite[Thm.~3.2]{S24}. \hfill $\Box$

\medskip

Let $M=\H^1_\cris(A/W)$ be the contravariant Dieudonn\'e module of $A$.
Recall that $\bU_A=\bU_{\wedge^2M}$.

\subsubsection*{The non-supersingular case}

Let $n\geqslant m\geqslant 1$.
Precomposing a map $M^\vee/p^m\to M/p^m$ with the surjective map 
$M^\vee/p^{n}\to M^\vee/p^m$ and composing with the injective map
$M/p^m\to M/p^{n}$ given by multiplication by $p^{n-m}$, we obtain an embedding
$$\Hom_\E(M^\vee/p^m, M/p^m)\hookrightarrow 
\Hom_\E(M^\vee/p^{n}, M/p^{n}).$$
Taking the anti-self-dual parts we get an embedding
$$\Hom_\E(M^\vee/p^m, M/p^m)^\skew\hookrightarrow 
\Hom_\E(M^\vee/p^{m+1}, M/p^{m+1})^\skew.$$
By Corollary \ref{2.14} the $p$-exponent of $\bU_A$ 
is the least positive integer $m$ such that 
$$\dim\,\Hom_\E(M^\vee/p^m, M/p^m)^\skew=\dim\,\Hom_\E(M^\vee/p^{n}, M/p^{n})^\skew.$$ 
Thus $\pexp(A)\leqslant r$, where $r$ is the least positive integer such that for all $s\geqslant r$ we have
$$\dim\,\Hom_\E(M^\vee/p^r, M/p^r)^{\skew}
=\dim\,\Hom_\E(M^\vee/p^s, M/p^s)^\skew.$$ 

For any Dieudonn\'e modules $N$ and $N'$
there exists an integer $f_{N,N'}\geqslant 0$ such that the image of 
the reduction modulo $p^m$ map
$$\Hom_\E(N/p^n,N'/p^n)\to \Hom_\E(N/p^m,N'/p^m)$$
is finite if and only if $n\geqslant m+f_{N,N'}$, see \cite[\S 6.1]{GV}.
From \cite[\S 6.1 (i), (ii)]{GV}, see also \cite[Lemma 7.8]{LNV},
we obtain $r=f_{M^\vee,M}$.

The abelian varieties $A$ and $A^\vee$ have the same Newton polygon that we denote by 
$\nu$. By \cite[Corollary 9.2]{LNV} we have $f_{M^\vee,M}\leqslant 2\nu(g)$. If $A$
is not supersingular, we have $\nu(g)<g/2$, hence $f_{M^\vee,M}\leqslant g-1$.
This finishes the proof in the non-supersingular case.

\subsubsection*{The supersingular case} 

Let $A$ be a supersingular abelian variety over $k$ of dimension $g\geqslant 1$. 
Let $E$ be a supersingular elliptic curve over $k$ and let $M_E$ be its Dieudonn\'e module. 
We have $M_E\cong\E/(F-V)\E$.
Since $A$ is isogenous to $E^g$, we have an isomorphism of $\E$-modules 
$M_\Q\cong M_{E,\Q}^g$.
Recall that an abelian variety is called superspecial if it is isomorphic
to $E^n$ for some $n\geqslant 1$. We call a Dieudonn\'e module superspecial if it is
isomorphic to $M_E^n$ for some $n\geqslant 1$. Recall that the $a$-number of $A$ is defined as
$$
a(A):=\dim_k(\Hom_k(\alpha_p,A))=\dim_k(M/(FM+VM)),
$$
hence the inequality $a(M)\leqslant g$.

Let us recall the basic theory of supersingular Dieudonn\'e modules as
developed in \cite[\S 1]{Li89}. By \cite[Lemma 1.3]{Li89} there exists a smallest
superspecial Dieudonn\'e module in $M_\Q\cong M_{E,\Q}^g$ containing $M$.
Let us call it the {\em superspecial envelope} of $M$ and denote it by $S'(M)$, using
the notation of \cite[\S 1]{Li89}.
The $W$-module $S'(M)/M$ has finite length.

The $W$-module $S'(M)$ is the direct sum of $g$ copies of $M_E=We\oplus Wf$,
where $Fe=Ve=f$ and $Ff=Vf=pe$. 
Consider the filtration by the images of powers of Frobenius 
$F^i S'(M)\subset S'(M)$ for $i\geqslant 0$, where
$F^0 S'(M)=S'(M)$. 
Note that the same filtration is obtained by using $V$ instead of $F$.
Thus each of $F$ and $V$ induces a ($\si$-linear or $\si^{-1}$-linear) bijective homomorphism 
of $g$-dimensional $k$-vector spaces:
$$F^iS'(M)/F^{i+1}S'(M)\tilde\lra F^{i+1}S'(M)/F^{i+2}S'(M).$$
This implies the injectivity of each of $F$ and $V$ on the graded factors of
the induced filtration $M^i:=M\cap F^i S'(M)$:
$$ M^i/M^{i+1} \hookrightarrow M^{i+1}/M^{i+2}.$$
Let $n$ be the smallest integer such that $F^n S'(M)\subset M$, which exists because $F$ is topologically nilpotent. Equivalently, 
$n$ is the smallest integer such that
$M^n/M^{n+1}\simeq k^g$. (For all $i\geqslant n$ we have $M^i/M^{i+1}\simeq k^g$.)

Write $s_i=\dim_k(M^i/M^{i+1})$, so $s_n=g$.
Following \cite[\S 1]{Li89} we claim that if $s_i<g$, then 
$$F(M^i/M^{i+1})\neq V(M^i/M^{i+1}).$$
In this case $M^{i+1}/M^{i+2}$ contains two distinct subspaces of dimension $s_i$, hence 
we have $s_i< s_{i+1}$ for $i\leqslant n-1$.

Let us prove the claim. By \cite[Lemma 1.6]{Li89}, there is an element $v\in M$ 
whose $e_i$-coefficients reduced modulo $p$ are linearly independent over $\F_{p^2}$. 
(Equivalently, the image of $v$ in $M/M^1$ is not contained in $L\otimes_{\F_{p^2}}k$
for a $\F_{p^2}$-vector space $L$ of dimension at most $g-1$. If this is not the case,
then $M$ is contained in
a proper superspecial submodule of $S'(M)$, which contradicts the fact that $S'(M)$
is the superspecial envelope of $M$, see \cite{Li89} for details.) 
If $s_i<g$ and $F(M^i/M^{i+1})=V(M^i/M^{i+1})$, then 
this is a proper $k$-vector subspace of the $k$-vector space $M^{i+1}/M^{i+2}$
of dimension at most $g-1$,
which is invariant under $x\to x^{p^2}$, and so it is $L\otimes_{\F_{p^2}}k$
for some $\F_{p^2}$-vector space $L$. Then $F^{i+1}(v)\in L\otimes_{\F_{p^2}}k$,
but this contradicts the definition of~$v$.

Write $N=(F,V)M$ and consider the induced filtration $N^i:=N\cap F^i S'(M)$.
We have $(F,V)M\subset F S'(M)=VS'(M)$, so $N=N^1$.
Define $t_i=\dim_k(N^{i+1}/N^{i+2})$ for $i=0,\ldots,n-1$. The crucial inequality
is $t_i\geqslant s_i+1$ for $i=0,\ldots, n-1$.
This follows from the claim above. Indeed, for these values of $i$ we have two injective maps
$M^i/M^{i+1} \hookrightarrow M^{i+1}/M^{i+2}$ (given by $F$ and $V$)
whose images are different and contained in
$N^{i+1}/N^{i+2}\subset  M^{i+1}/M^{i+2}$.

We finally note that $F^n S'(M)\subset M$ implies $F^{n+1} S'(M)\subset (F,V)M$,
hence $N^{n+1}=F^{n+1} S'(M)=M^{n+1}$.

Since the $a$-number of $A$ is $\dim_k(M/(FM+VM))$, we obtain
$$a(A)={\rm length}(M/M^{n+1})-{\rm length}(N/N^{n+1})=
\sum_{i=0}^n s_i-\sum_{i=0}^{n-1} t_i=s_n-\sum_{i=0}^{n-1}(t_i-s_i)\leqslant g-n.$$

Let $A'_i$ be an abelian variety with Dieudonn\'e module $M^i$, for $i=0,\ldots,n$.
In particular, $A'_0=A$ and $A'_n\simeq E^g$ is superspecial. Thus we get
a chain of isogenies
$$A=A'_0\to A'_1\to \ldots \to A'_n\simeq E^{g},$$
such that the kernel $K_i$ of $A'_i\to A'_{i+1}$ is $(\alpha_p)^{s_i}$.
It is easy to see that $a(A'_{i+1})-a(A'_i)=t_i-s_i\geqslant 1$. It follows that $n\leqslant g-1$,
and $n=g-1$ if and only if $a(A)=1$ and $a(A'_i)=i+1$ for $i=0,\ldots,g-1$. 
In this last case we have $s_i=i+1$ for $i=0,\ldots,g-1$, so that $K_i=\alpha(A'_i)$
is the group subscheme of $A'_i$ generated by the images of morphisms $\alpha_p\to A'_i$.

Applying previous considerations to the dual abelian variety, we obtain a chain of isogenies
$$E^g\simeq A_n\to A_{n-1}\to \ldots \to A_0=A,$$
whose successive kernels are $\alpha$-groups. Again, we have $n\leqslant g-1$,
and $n=g-1$ if and only if $a(A)=1$ and $a(A_i)=i+1$ for $i=0,\ldots,g-1$. 
In this case, the kernel of $A_i\to A_{i-1}$ is $(\alpha_p)^i$, so in particular the kernel
of $E^g\to A_{g-2}$ is $(\alpha_p)^{g-1}$.

\medskip

\noindent{\em Proof of Theorem} \ref{ss}. 
By Proposition \ref{livia} and the fact that $p\bU_{E^g}=0$, 
to prove that the $p$-exponent of $\bU_A$ is at most $g-1$, 
we can assume that $n=g-1$, so that the kernel of $E^g\to A_{g-2}$
is $(\alpha_p)^{g-1}$.
Theorem \ref{ss} will follow if we show that
$\bU_{A_{g-2}}\to \bU_{E^g}$ is the zero map.

We start with the case $g=2$.

\ble \label{5.2}
Let $E^2\to A$ be an isogeny with kernel $(\alpha_p)^n$, where $n>0$.
Then the induced map $\bU_A\to \bU_{E^2}$ is zero.
\ele
{\em Proof.} 
There is a commutative diagram with exact rows
$$\xymatrix{\wedge^2(M_{E^2})\ar[r]&\bU_{E^2}(k)\ar[r]&0\\
\wedge^2(M_A)\ar[r]\ar[u]&\bU_A(k)\ar[u]\ar[r]&0}$$
The Dieudonn\'e module $M_{E^2}=(M_E)^{\oplus 2}$ has $W$-basis
$e_1, e_2, f_1, f_2$ such that $Fe_i=Ve_i=f_i$ and $Ff_i=Vf_i=pe_i$ for $i=1,2$.
Thus $\Nyg(\wedge^2(M_{E^2}))$ is generated by
$pe_1\wedge e_2$, $f_1\wedge f_2$, and $e_i\wedge f_j$ for $i,j\in\{1,2\}$.
Calculating the map 
$$F/p-1\colon \Nyg(\wedge^2(M_{E^2}))\to \wedge^2(M_{E^2})$$
shows that it is surjective onto $\Nyg(\wedge^2(M_{E^2}))$,
hence $\wedge^2(M_{E^2})\to \bU_{E^2}(k)$ can be identified
with the projection onto the 1-dimensional $k$-vector space 
$k e_1\wedge e_2 $. By commutativity of the diagram, it is enough to show that the projection of $\wedge^2(M_A)$ onto $ke_1\wedge ke_2$ is zero.

We have $n=2$ or $n=1$. In the first case, $M_A$ is the $W$-submodule of $M_{E^2}$
generated by $pe_1, pe_2, f_1, f_2$, hence the projection of $\wedge^2(M_A)$ onto $ke_1\wedge ke_2$ is zero. In the second case,
$M_A$ is the $W$-submodule of $M_{E^2}$ generated by $pe_1, pe_2, f_1, f_2$,
and an element $a_1e_1+a_2e_2$, where at least one of $a_1, a_2$ is non-zero modulo $p$. Once again the projection of $\wedge^2(M_A)$ onto $ke_1\wedge ke_2$ is zero. \hfill $\Box$

\bpr\label{a=g-1}
The map $\bU_{A_{g-2}}\to \bU_{E^g}$ is zero.
\epr
{\em Proof.} Write $E^g=\prod_{i=1}^g E^{(i)}$, where the $E^{(i)}\simeq E$ are isomorphic
supersingular elliptic curves. It is known that 
$$\bU_{E^g}=\bigoplus_{i<j} p_{ij}^*(\bU_{E^{(i)}\times E^{(j)}})\cong (\bG_{k,a})^{g(g-1)/2},$$
where $p_{ij}\colon E^g\to E^{(i)}\times E^{(j)}$ is the projection, see 
\cite[Prop.~4.14]{Yang}.
Thus it is enough to prove that the composition 
$\bU_{A_{g-2}}\to \bU_{E^g}\to \bU_{E^{(i)}\times E^{(j)}}$ is zero.
(The second map is induced by the natural morphism $E^{(i)}\times E^{(j)}\to E^g$.)

Let $B$ be the image of $E^{(i)}\times E^{(j)}$ in $A_{g-2}$. We have a commutative diagram
$$\xymatrix{E^g\ar[r]&A\\
E^{(i)}\times E^{(j)}\ar[u]\ar[r]&B\ar[u]}$$
It is enough to prove that the induced map $\bU_B\to \bU_{E^{(i)}\times E^{(j)}}$ is zero.
The kernel of $E^{(i)}\times E^{(j)}\to B$ is the intersection of the kernel of 
$E^g\to A_{g-2}$, which is isomorphic to $(\alpha_p)^{g-1}$, with $E^{(i)}\times E^{(j)}$.
If this kernel was zero, we would get an injective map $(\alpha_p)^{g-1}\hookrightarrow E^{g-2}$,
which is not possible. Thus this kernel is $\alpha_p$ or $(\alpha_p)^2$. 
It remains to apply Lemma \ref{5.2} which says that
$\bU_B\to\bU_{E^{(i)}\times E^{(j)}}$ is the zero map. \hfill $\Box$

\subsection{The $p$-exponent of $\bU_A$ can be $g-1$} \label{exa}

The aim of this section is to prove the following statement.

\bpr \label{3.6}
Assume that $p\geqslant 3$.
For every $g\geqslant 1$ there exists a 
supersingular abelian variety 
$A$ of dimension $g$ over $k$ such that $\bU_A$ is isogenous to $\prod_{i=1}^{g-1}\bW_i$.
In particular, the $p$-exponent of $\bU_A$ is $g-1$.
\epr
{\em Proof.} By Corollary \ref{2.14} it is enough to construct a Dieudonn\'e module $M$ 
isoclinic of slope $1/2$, with a principal 
quasi-polarisation, and such that the connected component of the quasi-algebraic group
$\End_\E(M/p^n)^\skew$ is isogenous to $\prod_{i=1}^{g-1} \bW_i$, for all large enough $n$.
Here the superscript `$\skew$' denotes the subgroup
of anti-self-dual endomorphisms with respect to the principal quasi-polarisation.
Indeed, let $S'(M)$ be the superspecial envelope of $M$.
We have $S'(M)\cong M_{E^g}$, where $E$ is a supersingular elliptic curve over $k$. 
Dualising the embedding $M\to S'(M)$, we obtain an isogeny of $p$-divisible
groups $E^g[p^\infty]\to D$. Define $A$ as the quotient of $E^g$ by the kernel of this isogeny.
Then $\dim(A)=g$ and the Dieudonn\'e module of $A$ is $M$. 

Following \cite[Example 3.3]{NV07}, we consider $M=\E/(F^g-V^g)\E$. 
Thus $M$ is a free $W$-module of rank $2g$ with basis 
$$e_i=F^{i-1}, \ \  i=1,\ldots, g, \quad e_j=V^{2g+1-j},\ \  j=g+1,\ldots,2g.$$
We have $Fe_i=e_{i+1}$ for $i=1,\ldots, g$, and $Fe_i=pe_{i+1}$ for $i=g+1,\ldots, 2g$
(with the convention that subscripts are considered as integers modulo $2g$.) 
By \cite[Lemma 5.8]{LNV}, the Dieudonn\'e module $M$ is isoclinic of slope $1/2$.

The Dieudonn\'e module $M$ has a principal quasi-polarisation,
that is, a perfect, alternating, bilinear form
$$\psi\colon M\otimes_W M\to W,$$ 
such that $\psi(F(x),y)=\psi(x,V(y))^\si$,
defined as follows. Let $\theta\in W$
be an invertible element such that $\si^g(\theta)=-\theta$. Then
$$\psi(e_i,e_j)=0, \ \ i-j\neq g, \quad \psi(e_i,e_{i+g})=\si^{i-1}(\theta), \ \ i=1,\ldots, g,$$
is a principal quasi-polarisation. Write $B$ for the matrix of 
$\psi$ in the basis $e_1,\ldots,e_{2g}$, so that $B_{ij}=\psi(e_i,e_j)$ for $i,j=1,\ldots,2g$.

Since $M$ is generated by $e_1$ as a $\E$-module, there is at most one endomorphism
$f\in\End_\E(M/p^n)$ such that 
$$f(e_1)=\sum_{i=1}^{2g} a_i e_i,$$
for given $a_i\in W_n=W/p^n$ for $i=1,\ldots,2g$, and there is exactly one such $f$
if and only if $F^g f(e_1)=V^g f(e_1)$. We have
$$F^g f(e_1)=\sum_{i=1}^g a_i^{\si^g}p^{i-1}e_{i+g} +
\sum_{i=g+1}^{2g}a_i^{\si^g}p^{2g+1-i}e_{i-g},$$
$$V^g f(e_1)=\sum_{i=1}^g a_i^{\si^{-g}}p^{i-1}e_{i+g} +
\sum_{i=g+1}^{2g}a_i^{\si^{-g}}p^{2g+1-i}e_{i-g}.$$
Thus $F^g f(e_1)=V^g f(e_1)$ is equivalent to 
the condition that the image of $a_i$ in $W_{n+1-i}$ is invariant under $\si^{2g}$,
for $i=1,\ldots,2g$. Thus the connected component of $\End_\E(M)$ is isomorphic to 
\begin{equation}
\prod_{i=1}^g \bW_{i-1}\times \prod_{i=g+1}^{2g}\bW_{2g+1-i},
\label{W}\end{equation}
where the $i$-th factor is a group subscheme of $\bW_n$ with coordinate $a_i$.

Let $A$ be the $(2g\times 2g)$-matrix of $f$ in the basis $\{e_1,\ldots, e_{2g}\}$,
that is, $A_{ij}$ is the $e_i$-th coordinate of $f(e_j)$.
To $f$ we associate the map $\tilde f\in\Hom(M/p^n,(M/p^n)^\vee)$ that sends
$x$ to the linear form on $M/p^n$ whose value on $y$ is $\psi(Ax,y)$.
The dual map $\tilde f^\vee\in\Hom(M/p^n,(M/p^n)^\vee)$ sends $x$ to the linear form
whose value on $y$ is $\psi(Ay,x)=-\psi(x,Ay)$. Thus $\tilde f^\vee=-f^\vee$ if and only if
$A^t=BAB^{-1}$, or equivalently, $A=B^{-1}A^tB$. Since $f$ is an endomorphism of $\E$-modules, this condition is equivalent to 
$A(e_1)=B^{-1}A^tB(e_1)$. We have $B(e_1)=-\theta e_{g+1}$. Next, 
we use that $F^i(e_1)=e_{1+i}$ for $i=1,\ldots,g$ and $V^j(e_1)=e_{2g+1-j}$ 
for $j=1,\ldots, g$ to calculate
$$A^t(e_{g+1})=a_{g+1}e_1+a_g^\si e_2+\ldots+a_1^{\si^g}e_{g+1}+
a_{2g}^{\si^{-(g-1)}} e_{g+2}+\ldots+a_{g+2}^{\si^{-1}}e_{2g}.$$
From this we obtain 
$$B^{-1}A^tB(e_1)=a_1^{\si^g}e_1+
\sum_{i=2}^g \frac{\theta}{\theta^{\si^{i-1}}}a_{2g+2-i}^{\si^{i-g-1}}e_i
-\sum_{i=g+1}^{2g}\frac{\theta}{\theta^{\si^{i-g-1}}}a_{2g+2-i}^{\si^{i-g-1}}e_i.$$
Thus the involution on $\End_\E(M/p^n)$ that sends $A$ to 
$B^{-1}A^tB$ (which is indeed an involution because $B^t=-B^{-1}$)
sends the $i$-th factor of (\ref{W}) isomorphically onto the $(2g-i)$-th factor, and acts on
the $g+1$-th factor $\bW_g$ as $-1$. Since $p\neq 2$, we conclude that the connected
component of $\End_\E(M/p^n)$ is isomorphic to $\prod_{i=1}^{g-1} \bW_i$. \hfill $\Box$

\subsection{Abelian threefolds} \label{3-folds}

In this section we classify the isogeny classes of $\bU_A$ for abelian threefolds $A$
over an algebraically closed field $k$ of characteristic $p\geqslant 3$. 
 Recall that if $A$ is ordinary, then the slopes of the Newton polygon of $A$
are $(0,0,0,1,1,1)$, and if $A$ is supersingular, then all slopes are equal to $1/2$. There 
are three more cases, which we call {\em almost ordinary} if the slopes are $0,0,1/2,1/2,1,1$, 
{\em almost supersingular} if the slopes are $0,1/2,1/2,1/2,1/2,1$, and {\em $1/3$-type}
if the slopes are $1/3,1/3,1/3,2/3,2/3,2/3$. 

Let $\mathrm{BT}_1$ denote the category of $1$-truncated Barsotti--Tate groups, that is, 
finite commutative group $k$-schemes annihilated by $p$ such that
$\mathrm{Ker}(F)=\mathrm{Im}(V)$ and $\mathrm{Ker}(V)=\mathrm{Im}(F)$. 
A $1$-truncated Barsotti--Tate group is said to be \emph{indecomposable} if it is not
a direct product of two non-trivial subgroup schemes.

Let us recall Kraft's classification of $\mathrm{BT}_1$-group $k$-schemes \cite{K75}, see
\cite[\S 5.1]{Oo3}.
To every \emph{circular word} $w$ composed of letters $F$ and $V$ one
associates a $\mathrm{BT}_1$-group $G_w$ given by its covariant Dieudonné module $M_w$
which is an $\E$-module of finite length annihilated by $p$, in particular it is a 
finite dimensional $k$-vector space.
For example, the circular word $w=[FVF^2V^2]$ gives rise to the Dieudonn\'e module $M_w$
with a $k$-basis $\{e_1,e_2,e_3,e_4,e_5,e_6\}$ where $F$ and $V$ act as follows:
\begin{equation}
e_1\xrightarrow{F}e_2 \xleftarrow{V} e_3 \xrightarrow{F}e_4 \xrightarrow{F} e_5
\xleftarrow{V} e_6\xleftarrow{V} e_1.\label{eins}
\end{equation}
This means that
\begin{equation}\begin{aligned}
Fe_1=e_2, \ Ve_3=e_2,\ Fe_3=e_4, \ Fe_4=e_5, \ Ve_6=e_5, \ Ve_1=e_6;\\
Fe_2=Ve_2=Ve_4=Fe_5=Ve_5=Fe_6=0. \label{zwei}
\end{aligned}\end{equation}
The concatenation of two or more words corresponds to the direct sum of Dieudonné modules.
We say that a word $w$ is {\em indecomposable} if it is not of the form $w=w'^d$ for some word $w'$ and $d>1$. 

\bthe {\bf (Kraft)}
For every indecomposable circular word $w$, the $\mathrm{BT}_1$-group scheme $G_w$ is indecomposable. 
Every indecomposable $\mathrm{BT}_1$-group scheme is of the form $G_w$ for some indecomposable circular word $w$. 
\ethe
{\em Proof.} See, for example, \cite[Theorem 5.3]{Oo3}. \hfill $\Box$

\bpr\label{supergeneralAp}
Let $A$ be an abelian variety over $k$ of dimension $g\geqslant 1$ with $a$-number $1$
and $p$-rank $0$. Then the following properties hold.

{\rm (i)} There is an isomorphism of group $k$-schemes
$A[p]\cong G_w$, where $w=[V^gF^g]$.

{\rm (ii)} Assume, moreover, that $p\geqslant 3$. Then we have $\dim(\bU_A[p])=g-1$. 
\epr
{\em Proof.} (i) Let $M$ be the covariant Dieudonné module of $A[p^\infty]$.
 Since the $a$-number of $A$ is 1, we have $\dim_k(M/(FM+VM))=1$.
This implies that the $\E$-module $M$ is generated by one element.
Since the $p$-rank of $A$ is $0$, both $F$ and $V$ are nilpotent on $M/pM$.
Thus we can apply \cite[Lemma 5.7]{LNV}, which says that 
there is an isomorphism of $\E$-modules $M\cong\E/f\E$, where
\[f=a_0F^g+a_1F^{g-1}+\ldots+a_g+b_1V+\ldots+b_{g-1}V^{g-1}+b_gV^g\in \E\]
with $a_0,b_g\in W^\times$ and all other coefficients $a_i, b_j$ in $pW$. 
Reducing modulo $p$, we obtain an $\E$-module isomorphic to
the Dieudonné module $M_w$ for $w=[V^gF^g]$. 

(ii) By part (i), the $\mathrm{BT}_1$-group $k$-schemes 
$A[p]$ of abelian varieties $A$ of dimension $g$ with $a$-number $1$
and $p$-rank $0$ are isomorphic. By Corollary \ref{2.14} the dimension of 
$\bU_A[p]$ does not depend on the choice of $A$, and so can be calculated 
for the abelian variety with Dieudonn\'e module $\E/(F^g-V^g)\E$. Now (ii) follows from
Proposition \ref{3.6}. \hfill $\Box$

\bpr \label{24july}
The following table lists all possibilities for the $\mathrm{BT}_1$-group $k$-scheme
$A[p]\cong G_w$, where $A$ is an abelian threefold over $k$:
    \begin{center}{\rm
    \begin{tabular}{ c|c|c } \label{Table2}
    Newton polygon &    $a$      & $w$           \\ 
\hline
    supersingular &    $3$      &  $[FV]^3$  \\ 
    supersingular &    $2$      & $[FV][F^2V^2]$,  $[FVF^2V^2]$,  $[VFV^2F^2]$  \\    
    supersingular &    $1$      & $[F^3V^3]$       \\            
    $1/3$ type    &    $2$      &  $[F^2V][V^2F]$    \\ 
    $1/3$ type   &    $1$      & $[F^3V^3]$       \\                 
    almost supersingular& $2$   & $[F][V][FV]^2$ \\ 
    almost supersingular & $1$   & $[F][V][F^2V^2]$   \\           
     almost ordinary   & $1$   & $[F]^2[V]^2[FV]$ \\ 
     ordinary      &    $0$   & $[F]^3[V]^3$\\           
    \end{tabular}}
\end{center}
\epr
{\em Proof.} The superspecial case ($a=3$), the ordinary case, 
the almost ordinary case, and the almost supersingular cases are clear. 
The supergeneral case and the $1/3$ type case with $a=1$ follow from Proposition \ref{supergeneralAp} (i).

For the $1/3$ type case with $a=2$ and the supersingular case with $a=2$, 
the $p$-rank is $0$. There are only four possibilities for the circular words that give rise
to such group schemes: 
\begin{equation}
[FVF^2V^2], \ \ [VFV^2F^2], \ \ [FV][F^2V^2], \ \ [F^2V][V^2F]. \label{words}
\end{equation}

The group scheme $G_w$ associated to $w=[F^2V][V^2F]$ 
is minimal in the sense of \cite[Definition, p.~1024]{O05}.
In fact, we have $G_w\cong H_{1,2}[p]\oplus H_{2,1}[p]$, 
where $H_{c,d}$ is the $p$-divisible group
with contravariant Dieudonn\'e module $M_{c,d}$ as in the proof of
Proposition \ref{5.8}. By minimality, $A[p]\cong G_w$ implies that 
$A[p^\infty]\cong H_{1,2}\oplus H_{2,1}$,
so that the Newton polygon of $A$ is of $1/3$ type. 

If $w$ is one of the three words in (\ref{words}) other than $[F^2V][V^2F]$, then
\begin{equation}
\dim(\Hom_\E(M_w,M_w^\vee)^{\skew})=3. \label{d4}
\end{equation} 
When $A[p]\cong G_w$, 
from Corollary \ref{2.14} we obtain that $\dim(\bU_A)\geqslant\dim(\bU_A[p])=3$. 
For abelian varieties $A$ of dimension $g$ we have 
$\dim(\bU_A)\leqslant g(g-1)/2$ by (\ref{dim}), with equality
 if and only if $g=1$ or $A$ is supersingular (Remark \ref{ferry} (i)). Thus (\ref{d4}) implies
that our abelian threefold $A$ is supersingular.

The third named author proved (\ref{d4}) in his thesis, see \cite[Section 4.4]{YuanThesis}. 
For the convenience of the reader let us give this calculation here.

Let $w=[FVF^2V^2]$. Let 
$e_1^*,\ldots, e_6^*$ be the dual basis of the dual Dieudonn\'e module $M_w^\vee$. Dualising (\ref{eins}) we obtain
\begin{equation}
e_1^*\xleftarrow{V}e_2^* \xrightarrow{F} e_3^* \xleftarrow{V}e_4^*   \xleftarrow{V}e_5^*
\xrightarrow{F}  e_6^*\xrightarrow{F}  e_1^*.\label{drei}
\end{equation}
A linear map of $k$-vector spaces $f\colon M_w \to M_w^\vee$ is a map of Dieudonn\'e modules
if and only if it commutes with $F$ and $V$. Write $f(e_j)=\sum a_{ij} e_i^*$ with $a_{ij}\in k$
such that $a_{ij}=-a_{ji}$ for all $i$ and $j$, in particular, $a_{ii}=0$.
We consider subscripts as integers modulo 6.

Commuting with $F$ is equivalent to the conjunction of 

(1) $a_{i+1,j+1}= a_{ij}^\si$, when $Fe_j=e_{j+1}$ and $Ve_{i+1}=e_i$;

(2) $a_{ij}=0$, when $Fe_j=0$ and $Ve_{i+1}=e_i$.

Commuting with $V$ is equivalent to the conjunction of 

(3) $a_{i+1,j+1}=a_{i,j}^\si$, when $Ve_{j+1}=e_j$ and $Fe_i=e_{i+1}$;

(4) $a_{ij}=0$, when $Ve_j=0$ and $Fe_{i-1}=e_i$.

\noindent In our case $(i,j)$ satisfies condition (1) for $i=2,5,6$ and $j=1,3,4$,
condition (2) for $i,j\in\{2,5,6\}$, condition (3) for $i=1,3,4$ and $j=2,5,6$, and 
condition (4) for $i,j\in\{2,4,5\}$. The set of pairs $(i,j)$ with $i\neq j$ is the disjoint union
of chains where we have an arrow $(i,j)\mapsto (i+1,j+1)$ if $a_{i+1,j+1}= a_{ij}^\si$.
Counting the chains that do not contain a pair $(i,j)$ satisfying (2) or (4), we see that there
are exactly six such chains. Taking skew-symmetry of $(a_{ij})$ into account, we deduce that
$\dim_k(\Hom(M_w,M_w^\vee)^{\skew})=3$.

The case $w=[VFV^2F^2]$ is entirely similar.

It remains to deal with the case $w=[FV][F^2V^2]$. Let $M$ and $N$ be the Dieudonn\'e modules
given by $[FV]$ and $[F^2V^2]$, respectively.
Thus $M_w=M\oplus N$, and we have
$$\Hom_\E(M_w,M_w^\vee)^{\skew}\cong \Hom_\E(M,M^\vee)^{\skew}\oplus
\Hom_\E(N,N^\vee)^{\skew}\oplus\Hom_\E(M,N^\vee).$$
Since $\bU_{\wedge^2(M)}=0$, we have $\Hom_\E(M,M^\vee)^{\skew}=0$.
Next, $\bU_{\wedge^2(N)}\cong\bG_{a,k}$ implies that $\Hom_\E(N,N^\vee)^{\skew}$ 
has dimension 1. Noting that $N$ is self-dual,
it remains to check that the dimension of $\Hom_\E(M,N)$ is 2.
Since $M=\E/(F-V)\E$, this dimension is equal to the dimension of $N^{F=V}$. 
The Dieudonn\'e module $N$
is a $k$-vector space with basis $e_1,e_2,e_3,e_4$ with action of $F$ and $V$ given by
$$e_1\xrightarrow{F}e_2 \xrightarrow{F} e_3 \xleftarrow{V}e_4 \xleftarrow{V} e_1.$$
Write $u=a_1e_1+a_2e_2+a_3e_3+a_4e_4\in N$. Then
$F(u)=a_1^\si e_2+a_2^\si e_3$ and $V(u)=a_1^{\si^{-1}}e_4+a_4^{\si^{-1}}e_3$. 
Hence $N^{F=V}$ is given by $a_1=0$ and $a_4=a_2^{\si^2}$. Thus $a_2$ and $a_3$
are free variables, so the dimension of $N^{F=V}$ is 2. This finishes the proof. \hfill $\Box$

\bpr\label{NPPabe3table}
Let $A$ be an abelian threefold over an algebraically closed field $k$ of characteristic $p\geqslant 3$. The following table lists all possible isogeny classes of $\bU_A$:
    \begin{center}{\rm
    \begin{tabular}{ c|c|c|c|c } \label{Table2}
    Newton polygon &    $a$      & $\dim(\bU_A[p])$& $\dim(\bU_A)$ & $\bU_A$ up to isogeny\\ 
\hline
    supersingular &    $3$      &  $3$       & $3$      & $(\bG_{a,k})^3$  \\ 
    supersingular &    $2$      &  $3$       & $3$     & $(\bG_{a,k})^3$\\    
    supersingular &    $1$      &  $2$       & $3$     & $\bG_{a,k}\times \bW_2$\\            
    $1/3$ type    &    $2$      &  $2$       & $2$     & $(\bG_{a,k})^2$  \\ 
    $1/3$ type    &    $1$      &  $2$       & $2$     & $(\bG_{a,k})^2$\\                 
    almost supersingular& $2$   &  $1$       & $1$     & $\bG_{a,k}$  \\ 
    almost supersingular & $1$   &  $1$       & $1$     & $\bG_{a,k}$\\           
     almost ordinary   & $1$     &  $0$       & $0$     & $0$  \\ 
     ordinary      &    $0$      &  $0$       & $0$     & $0$\\           
    \end{tabular}}
\end{center}
In particular, the isogeny class of $\bU_A$ is uniquely determined by the Newton polygon and the
$a$-number of $A$. We have
$\bU_A\simeq(\bG_{a,k})^n$, where $0\leqslant n\leqslant 3$, unless $A$ is supergeneral, in which case
$\bU_A$ is isogenous to $\bG_{a,k}\times \bW_2$. 
\epr
{\em Proof.} The superspecial case is clear. Using the dimension formula (\ref{dim}),
the supergeneral case and the $1/3$ type case 
with $a=1$ follow from Proposition \ref{supergeneralAp} (ii). 
We dealt with the supersingular case with $a=2$ in the proof of Proposition \ref{24july}
which shows that $pU_A=0$ in this case. 

For the $1/3$ type case with $a=2$, we have $A[p]\cong G_w$ where $w=[F^2V][V^2F]$.
As was observed in the proof of Proposition \ref{24july}, this implies that
$A[p^\infty]\cong H_{1,2}\oplus H_{2,1}$. In this case we have $pU_A=0$,
as follows from the proofs of Proposition \ref{5.8} and Corollary \ref{5.10}.

The remaining cases were handled in Remark \ref{ferry} (ii) and (iii).
\hfill $\Box$

\section{Ekedahl--Oort types and $\dim(\bU_A[p])$} \label{EO}

Let $A$ be an abelian variety of dimension $g$ over an algebraically closed field $k$
of characteristic $p\geqslant 3$.
In this section we compute $\dim(\bU_A[p])$
when $A$ is principally polarised in terms of the Ekedahl--Oort type of $A$.

In \cite{Oo1}, Ekedahl and Oort stratified the moduli space of 
principally polarised abelian varieties of dimension $g$ over $k$ 
according to the isomorphism class of $A[p]$.
In total there are $2^g$ such isomorphism classes; they are numbered by the Ekedahl--Oort types. 
An Ekedahl--Oort type is a non-decreasing function 
$$\varphi\colon\{0,1,\ldots,g\}\to\{0,1,\ldots,g\}, \quad \varphi(0)=0,$$ such that 
$\varphi(i)-\varphi(i-1)$ is 0 or 1. Given a  subset $P\subset \{1,\ldots,g\}$
we define
$\varphi_P(i)=\varphi_P(i-1)+1$ if $i\in P$, and $\varphi_P(i)=\varphi_P(i-1)$ otherwise.
It is clear that every Ekedahl--Oort type can be written as $\varphi_P$ for some subset 
$P\subset \{1,\ldots,g\}$. 
The correspondence between Ekedahl--Oort types and isomorphism classes of $A[p]$ is recalled in Theorem \ref{oort}.
Note that the $a$-number of $A$ is $g-|P|$.

\bthe \label{E-O}
Let $A$ be a principally polarised abelian variety of dimension $g$ over an algebraically closed
field of characteristic $p\geqslant 3$.
If the Ekedahl--Oort type of $A$ is $\varphi_P$, where $P=\{m_1<\ldots<m_{g-a}\}$, 
where $a$ is the $a$-number of $A$, then 
    \[\dim (\bU_A[p])=\frac{a(a-1)}{2}+\sum_{i=1}^{g-a}(m_i-i).\]
\ethe

When $P=\emptyset$, the sum is understood to be 0.

\bexa{\rm \label{ex1} Let $A$ be a principally polarised abelian variety.

(1) Suppose that $A$ is supersingular.
Then we have $\dim(\bU_A)=g(g-1)/2$.
There are two extreme cases: 
\begin{itemize}
\item 
$A$ is superspecial if and only if $a=g$, which is equivalent to 
the condition that $\varphi$ is the zero function (that is, $P=\emptyset$).
Theorem \ref{E-O} gives $\dim(\bU_A[p])=g(g-1)/2$, so that $p\bU_A=0$. This is already proved in \cite[Proposition 4.4]{Yang}.
\item
$A$ is supergeneral if and only if $a=1$.
Then $P=\{2,\ldots,g\}$ and
$\varphi(i)=i-1$ for $i\geqslant 1$. (This can be deduced from \cite[(9.7)]{Oo1} using that 
the $p$-rank is $0$.) Theorem \ref{E-O} says that $\dim(\bU_A[p])=g-1$.
\end{itemize}

(2) Suppose that $A$ is ordinary. This is equivalent to $a=0$, and thus to $\varphi(i)=i$
for all $i=0,1,\ldots, g$.
}
\eexa

We now recall the construction of {\em covariant} self-dual Dieudonn\'e modules
associated to Ekedahl--Oort types  \cite[\S 9]{Oo1}. 
Let $P=\{m_1<...<m_h\}$ be a subset of $\{1,2,...,g\}$ of size $h=g-a$, and let 
$\varphi=\varphi_P$.
Extend $\varphi$ to a non-decreasing function 
$$\psi\colon\{0,1,\ldots,2g\}\to\{0,1,\ldots,g\}, \quad \psi(0)=0, \ \psi(2g)=g,$$ by setting
$\psi(i)=\varphi(2g-i)+i-g$ for $i=g+1,...,2g$. 
Let $m_1<\ldots<m_g$ be the increasing sequence of all integers 
$1\leqslant i\leqslant 2g$ such that $\psi(i)=\psi(i-1)+1$. Likewise, let
$n_1>\ldots>n_g$ be the decreasing sequence of all integers $1\leqslant j\leqslant 2g$
such that $\psi(j)=\psi(j-1)$. In particular, 
$\{1,...,g\}$ is the disjoint union of $P=\{m_1,\ldots,m_h\}$
and its complement $Q=\{n_g,\ldots,n_{h+1}\}$.

Let $M_\varphi$ be a vector space over $\F_p$ with basis $Z_1,\ldots,Z_{2g}$.
Let us write $\sZ_1:={\rm Span}(Z_1,\ldots,Z_g)$ and 
$\sZ_2:={\rm Span}(Z_{g+1},\ldots,Z_{2g})$.
We rename the set of basis vectors 
$$Z_{m_i}=X_i, \quad Z_{n_i}=Y_i, \ \ \text{where} \ \ i=1,\ldots,g,$$
and write $\sX:={\rm Span}(X_i)$, $\sY:={\rm Span}(Y_j)$.
The $\F_p$-vector space $M_\varphi$ has an alternating bilinear form $\langle x,y\rangle$
which is zero on $\sX$ and $\sY$,
and satisfies $\langle X_i,Y_j\rangle=\delta_{ij}$, where $\delta_{ij}$ is Kronecker's delta.
We have $m_i+n_i=2g+1$ for all $i=1,\ldots,g$. Using this one checks that 
$\sZ_1$ and $\sZ_2$ are maximal hyperbolic subspaces, alongside with $\sX$ and $\sY$.

Define $\mathcal F\in\End_{\F_p}(M_\varphi)$ as follows:
$$\mathcal F(X_i)=Z_i, \quad \mathcal F(Y_i)=0. $$
We have $\langle \mathcal Fx,\mathcal Fy\rangle=0$ for all $x,y\in M_\varphi$, because 
$\Im(\mathcal F)=\sZ_1$ which is maximal hyperbolic. 
Define  $\mathcal V\in\End_{\F_p}(M_\varphi)$ so that
\begin{equation}
\langle\mathcal F x,y\rangle=\langle x, \mathcal Vy\rangle \label{Y5}
\end{equation}
for all $x,y\in M_\varphi$. By the non-degeneracy of the form we have $\mathcal V\mathcal F=0$.
Using this non-degeneracy 
and the fact that $\Ker(\mathcal F)=\sY$, we see that 
$\Im(\mathcal V)=\sY$ and $\Ker(\mathcal V)=\sZ_1$. We deduce that
$\langle \mathcal Vx,\mathcal Vy\rangle=0$ for all $x,y\in M_\varphi$.
This implies $\mathcal F\mathcal V=0$.

\bthe\emph{\textbf{(Oort)}} \label{oort}
Let $k$ be an algebraically closed field of positive characteristic $p\neq 2$.
Let $A$ be a principally polarised abelian variety over $k$ of Ekedahl--Oort type $\varphi$. 
Then the covariant Dieudonné module $M$ of $A[p]$ is $M\cong M_\varphi\otimes k$ with 
the natural $\si$-linear action of $\mathcal F$ and the natural 
$\si^{-1}$-linear action of $\mathcal V$. 
The Dieudonn\'e module $M$ is self-dual with respect to the alternating
bilinear form $\langle x,y\rangle$.
\ethe
{\em Proof.} This is \cite[Theorem 9.4]{Oo1}. $\square$

\medskip

Let $E_{i,j}$ be the $(g\times g)$-matrix with $(i,j)$-entry 1 and all other entries 0.
Define
$$U=\sum_{i=1}^h E_{i,m_i}, \quad V=\sum_{j=h+1}^g E_{j,n_j}.$$
The matrix $U$ encodes that $\mathcal F(X_{m_i})=Z_{m_i}=X_i$, whereas
$V$ encodes that $\mathcal F(X_{n_i})=Z_{n_i}=Y_i$. Thus
$U+V$ is the permutation matrix of the permutation of $\{1,\ldots,g\}$
defined on $P$ as the unique order-preserving
function $P\to\{1,\ldots,h\}$ (the inverse of $i\mapsto m_i$)
and on $Q$ as the unique order-inversing function 
$Q\to \{h+1,\ldots,g\}$ (the inverse of $j\mapsto n_j$).

\ble
The group of self-dual endomorphisms of the self-dual Dieudonn\'e module
$M$ is isomorphic to the additive group of 
pairs of $(g\times g)$-matrices $(A,C)$ with entries in $k$ such that 
\begin{equation}
C^t=-C, \quad VA=C^{(p)}U+A^{t(p)} V, \quad UA=A^{(p)}U, \label{Y6}
\end{equation}
where superscript $(p)$ means that every entry of the matrix is raised to the power~$p$.
\ele
{\em Proof.}
We use the basis $\mathcal B=\{X_1,\ldots,X_g,Y_1,\ldots,Y_g\}$.
In this basis the matrix of the form $\langle x,y\rangle$ is 
$$\left(\begin{matrix}
    0&I\\
    -I&0
\end{matrix}\right).$$
The linear action of $\mathcal F$ and $\mathcal V$ on $M_\varphi$ in the basis 
$\mathcal B$ is given by the matrices
\[\mathcal [\mathcal F]=\left(\begin{matrix}
    U&0\\
    V&0
\end{matrix}\right),\quad\quad [\mathcal V]= \left(\begin{matrix}
    0&0\\
    -V^t&U^t
\end{matrix}\right).\]
Here $[\mathcal V]$ is obtained from $[\mathcal F]$ using (\ref{Y5}).

Let $f$ be an endomorphism of the self-dual Dieudonn\'e module $M$. 
We have $f=f^\vee$ if and only if the matrix of $f$ in the basis $\mathcal B$ is
\[[f]=\left(\begin{matrix}
    A&B\\
    C&D
\end{matrix}\right)\]
where $D=A^t$, $B^t=-B$, $C=-C^t$. The condition $\mathcal F f=f^\si \mathcal F$
is equivalent to
$$UA=A^{(p)}U+B^{(p)}V, \quad VA=C^{(p)}U+A^{t(p)}V, \quad UB=VB=0,$$
Since $U+V$ is invertible, we have $B=0$. This implies the lemma. $\Box$

\medskip

As an example, suppose that $P=\{1,\ldots,g\}$ so that $Q=\emptyset$. Then $U=I$ and $V=0$.
It follows that $C=0$ and $A$ is an arbitrary element of ${\rm Mat}_g(\F_p)$. This is the case
for ordinary abelian varieties. Now suppose that $P=\emptyset$. Then $U=0$ and $V$ is the 
matrix of the permutation $i\mapsto g+1-i$. Then $C$ is an arbitrary skew-symmetric
matrix with entries in $k$. Since $p\neq 2$, this gives a vector space over $k$ of dimension 
$g(g-1)/2$. The condition on $A$ gives a finite-dimensional vector space over $\F_p$.
This is the case for superspecial abelian varieties. Both results agree with Theorem~\ref{E-O}.

\medskip

\noindent{\em Proof of Theorem} \ref{E-O}. 
By Corollary \ref{2.14} it is enough to prove the following statement.

\medskip

\noindent ($\star$) 
{\em 
The subgroup $S\subset{\rm Mat}_g(k)$, defined by the property that for every $A\in S$ there is a 
matrix $C\in {\rm Mat}_g(k)$ with $(A,C)$ a solution of (\ref{Y6}), is isomorphic to a direct sum of a 
finite-dimensional $\F_p$-vector space and a $k$-vector space of dimension $\sum_{i=1}^h(m_i-i)$.}

\medskip

Indeed, by ($\star$), for every $A\in S$ there is a matrix $C$ such that $(A,C)$ satisfies (\ref{Y6}).
If $(A,C+C')$ is a solution of (\ref{Y6}), then
$C'U=0$ which says that the $h$ columns of the skew-symmetric matrix 
$C'$ corresponding to the elements of $P$ are zero. 
Thus $C$ is uniquely determined 
up to addition of an element of a $k$-vector subspace of ${\rm Mat}_g(k)$
of dimension $(g-h)(g-h-1)/2$, so ($\star$) implies the theorem. 
\medskip

\ble\label{solutions}
    The subgroup $S\subset{\rm Mat}_g(k)$ is given by the following conditions:

\smallskip

{\rm(1)} $a_{m_i,m_j}=a_{i,j}^p$ for $i,j=1,\ldots,h$;

{\rm(2)} $a_{n_j,n_i}=a_{i,j}^p$ for $i,j=h+1,\ldots,g$;

{\rm(3)} if $i\geqslant h+1$, $j\leqslant h$, then $a_{i,j}=0$;

{\rm(4)} if $i\in P$, $j \in Q$, then $a_{i,j}=0$. 
\ele
{\em Proof.} The $i$-th row of $UA$
is the $m_i$-th row of $A$, for $i=1,\ldots,h$, and all other rows are zero.
The $m_j$-th column of $A^{(p)}U$ is the $j$-th column of $A^{(p)}$, for $j=1,\ldots, h$,
and all other columns are zero. 
This implies that $UA=A^{(p)}U$ is equivalent to the conjunction of conditions (1), (3), (4).

The $m_j$-th column of $C^{(p)}U$ is the $j$-th column of $C^{(p)}$, 
for $j=1,\ldots, h$, and all other columns are zero.
The $n_j$-th column of $A^{t(p)} V$ is the $j$-th row of $A^{(p)}$, for $j=h+1,\ldots,g$,
and all other columns are zero.
The $i$-th row of $VA$ is the $n_i$-th row of $A$, for $i=h+1,\ldots,g$, and all other rows are zero.
From this one sees that the conditions imposed on $A$ by the equation $VA=C^{(p)}U+A^{t(p)} V$ are equivalent to the conjunction of (3) and (2). $\square$

\medskip

To prove ($\star$) we consider a directed graph
$\mathcal{G}_{\varphi}$, whose vertices are pairs $(i,j)$,
for $i,j=1,\ldots,g$, and the directed edges encode the relations of Lemma \ref{solutions}. 
Namely, the edges of $\mathcal{G}_{\varphi}$ are of the following types (loops are allowed):

(m) for $i,j=1,\ldots,h$ there is an edge from $(i,j)$ to $(m_i,m_j)$; 

(n) for $i,j=h+1,\ldots,g$ there is an edge from $(i,j)$ to $(n_j,n_i)$.

\noindent 
We say that a vertex $(i,j)$ is {\em marked} if $(i,j)$ satisfies condition (3) or (4)
of Lemma \ref{solutions}, that is, if $j\leqslant h<i$ or $(i,j)\in P\times Q$.

\ble \label{last}
The following properties hold:

    {\rm(1)} if there is an edge from $(i_1,j_1)$ to $(i_2,j_2)$, then 
$i_1-j_1$ and $i_2-j_2$ have the same sign (so the connected components `do not cross
the diagonal');

    {\rm(2)} every vertex has at most one incoming edge and at most one outgoing edge;

    {\rm(3)} the connected components of the graph $\mathcal{G}_{\varphi}$ 
are either cycles $(i_1,j_1)\rightarrow\ldots \rightarrow (i_r,j_r)\rightarrow (i_1,j_1)$
or chains $(i_1,j_1)\rightarrow\ldots \rightarrow (i_s,j_s)$ (without repetition of vertices);

    {\rm(4)} for every `diagonal' vertex $(i,i)$ there is exactly one incoming and one outgoing edge;
the edge starting at a diagonal vertex ends in a diagonal vertex;

    {\rm(5)} the connected components of the set of diagonal vertices are cycles;

    {\rm(6)} no edge originates at $(i,j)$ if and only if 
    $i \leqslant h < j$ or $j \leqslant h < i$;

    {\rm(7)} there are exactly $h(g-h)$ chains above (respectively, below) the diagonal;

    {\rm(8)} no cycle contains a marked vertex;

    {\rm(9)} every chain below the diagonal ends with a marked vertex; 

    {\rm(10)} every marked vertex above the diagonal is the start point of a chain;

    {\rm(11)} the number of chains without marked vertices is $\sum_{i=1}^h(m_i-i)$. 
\ele
{\em Proof.} For a type (m) edge, $i<j$ implies $m_i<m_j$ because the function 
$i\mapsto m_i$ is order-preserving. For a type (n) edge, $i<j$ implies $n_j<n_i$
because the function $i\mapsto n_i$ is order-inversing. This proves (1).
Properties (2), (4), (6) are immediate from the construction. Property
(3) is a consequence of (2). Property (5) is a consequence of (4). A chain is uniquely
determined by its end point $(i,j)$, which must be such that $i \leqslant h < j$ or $j \leqslant h < i$
by (6). Moreover, any such pair $(i,j)$ is the end point of some chain (which may have
length 1). This implies (7).
For (8), we note that a marked vertex is either $j\leqslant h<i$, in which case there is no
outgoing edge, or $(i,j)\in P\times Q$, in which case there is no incoming edge.
Property (9) is consequence of (6). 
Every marked vertex above the diagonal is of the form $(i,j)$, where $(i,j)\in P\times Q$, so
there is no incoming edge, thus $(i,j)$ is the start point of a chain. This proves (10).

Let us prove (11). All marked points are either start or end points of chains, so
a chain has no marked points if and only if its start point and its end point are not marked.
This means that it starts at a vertex $(i,j)\in Q\times P$, and is above the diagonal.
Thus a chain has no marked points if and only if it starts at a vertex
$(i,j)\in Q\times P$ such that $j>i$. The number of such pairs $(i,j)$ is 
$$\sum_{i=1}^h |\{n\in Q| n<m_i\}|=\sum_{i=1}^h(m_i-i).$$
This finishes the proof of the lemma. $\Box$

\medskip

\noindent{\em End of proof of Theorem} \ref{E-O}. Let us prove ($\star$).
By Lemma \ref{solutions} we need to compute the space of
$k$-valued functions $a_{ij}$ on the set of vertices of the directed graph $\mathcal{G}_{\varphi}$
such that if $(i,j)$ is a marked vertex then $a_{ij}=0$, and if there is an edge from
$(i_1,j_1)$ to $(i_2,j_2)$ then $a_{i_2,j_2}=a_{i_1,j_1}^p$. This space
is the direct sum of spaces given by the connected components of $\mathcal{G}_{\varphi}$.

The directed graph $\mathcal{G}_{\varphi}$ has the three kinds of connected components:
cycles (which contain no marked vertices), chains with marked vertices, and chains with
no marked vertices. The vector space associated to a cycle with $n$ vertices
is isomorphic to $\F_{p^n}$. The vector space associated to a chain with marked vertices is 0.
Finally, the vector space associated to a chain without marked vertices is isomorphic to $k$.
The proof finishes by an application of Lemma \ref{last} (11).
$\square$

\bco \label{d1}
Let $A$ be a principally polarised abelian variety over an algebraically closed
field $k$ of characteristic $p\geqslant 3$. Let $r=\dim (\bU_A[p])$.
Then there exists an abelian variety $B$ over $k$ of dimension at most $r+1$ such that
$\bU_A\cong \bU_B$. In particular, we have
$$\dim(\bU_A)\leqslant r(r+1)/2 \quad \text{and} \quad \pexp(\bU_A)\leqslant r.$$
\eco
{\em Proof.} Let $f$ be the $p$-rank of $A$. Since $k$ is perfect,
the $p$-divisible group $A[p^\infty]$ is isomorphic to the product of $p$-divisible
groups of slopes 0 and 1, and the $p$-divisible group
of an abelian variety $B$ of dimension $g-f$. 
By Corollaries \ref{dim 0} and \ref{d3} this implies that $\bU_A\cong \bU_B$.

By \cite[(9.7)]{Oo1}, $f$ is the maximal value of $i$ such that $\varphi(i)=i$, i.e.,
the sequence $\varphi(i)$ starts with the sequence $1,2,\ldots,f-1,f,f$.
Thus $f$ is the maximal value of $i$ such that $m_i=i$. In particular,
we have $m_i-i\geqslant 1$ for $i\geqslant f+1$. Now Theorem \ref{E-O} gives
$$r=\frac{a(a-1)}{2}+\sum_{i=1}^{g-a}(m_i-i)\geqslant \frac{a(a-1)}{2}+g-a-f=
\frac{a(a-3)}{2}+g-f\geqslant g-f-1,$$
where we used that $a(a-3)/2\geqslant -1$ for all $a\in\Z$. 
Now (\ref{dim}) gives  
$\dim(\bU_A)=\dim(\bU_B)\leqslant(g-f)(g-f-1)/2\leqslant r(r+1)/2$. \hfill $\Box$

\bexa \label{d2}
{\rm It is instructive to explore what Corollary \ref{d1} gives for small values of $r$.
Suppose that $\bU_A$ is isogenous to $\bW_{n_1}\times\ldots\times\bW_{n_r}$, where 
$n_1\leqslant\ldots\leqslant n_r$.

(1) If $r=1$, then $n_1=1$.

(2) If $r=2$, then $(n_1,n_2)$ is either $(1,1)$ or $(1,2)$.

(3) If $r=3$, then $n_1+n_2+n_3\leqslant 6$. There are two possible cases: $\dim(B)=3$
and $\dim(B)=4$. In the first case we have $(n_1,n_2,n_3)=(1,1,1)$. In the second case,
we have $n_i\leqslant 3$ by Theorem \ref{ss}. Listing all possibilities for the Newton polygon of $B$, we see that  $(n_1,n_2,n_3)=(1,1,2)$
is the only possibility in the non-supersingular case besides $(1,1,1)$.
In the supersingular case we have $n_1+n_2+n_3=6$.
In view of Proposition \ref{3.6} we expect that the case $(n_1,n_2,n_3)=(1,2,3)$ does occur,
but we are not sure if $(n_1,n_2,n_3)=(2,2,2)$ does.
}\eexa

\bco \label{4.7}
The following table gives all possible isogeny classes of $\bU_A$ for
principally polarised abelian threefolds over an algebraically closed field $k$ of characteristic
$p\geqslant 3$ of a given Ekedahl--Oort type:
{\rm \begin{center}
    \begin{tabular}{ c|c|c|c|c } \label{Table2}
     E-O type   $\varphi$  &$a$&$p$-rank  & Newton polygon & $\bU_A$ up to isogeny\\ 
\hline
     $0,0,0$   &3  &0     &   supersingular  &   $(\bG_{a,k})^3$  \\ 
     $0,0,1$   &2  &0     &   supersingular  &   $(\bG_{a,k})^3$\\       
     $0,1,2$   &1  &0     &   supersingular  &  $\bG_{a,k}\times \bW_2$\\         
     $0,1,2$   &1  &0&   $1/3$-type       &   $(\bG_{a,k})^2$\\              
    $0,1,1$   &2  &0    &   $1/3$-type     & $(\bG_{a,k})^2$  \\    
     $1,1,1$   &2   &1    &   almost supersingular    &  $\bG_{a,k}$  \\ 
     $1,1,2$  &1    &1    &   almost supersingular  &   $\bG_{a,k}$ \\           
     $1,2,2$   &1   &2    &   almost ordinary  &   $0$  \\ 
     $1,2,3$    &0   &3   &   ordinary  &   $0$\\           
    \end{tabular}
\end{center}}
\eco
{\em Proof.} The Ekedahl--Oort type of $\bU_A[p]$ determines $a$, the $p$-rank, and $\dim(\bU_A[p])$,
whereas the Newton polygon determines $\dim(\bU_A)$. 

Ekedahl--Oort types $(0,0,0)$ and $(1,2,3)$ were dealt with in Example \ref{ex1}.

If $a=2$, then $\varphi$ can be $(0,0,1)$, $(0,1,1)$, and $(1,1,1)$. 
A result of Oort \cite[Theorem (8.3)]{Oo1} says that in the first case $A$ is supersingular, in the second case $A$ of the $1/3$-type, and in the third case $A$ is almost supersingular.

If $a=1$, then $\varphi$ can be $(0,1,2)$, $(1,1,2)$, and $(1,2,2)$. We know from \cite[(9.7)]{Oo1}
that the $p$-rank $f$ of $A$ is the maximal value of $i$ such that $\varphi(i)=i$.
This implies that in the last case, the $p$-rank of $A$ is 2. Thus we have an isomorphism
of group $k$-schemes
$A[p]\cong(\Z/p)^2\times (\mu_p)^2\times G$, where $G$ is local-local.
It follows that $A$ is almost ordinary.
The same argument gives that in the second case the $p$-rank of $A$ is 1.
Then $A[p]\cong\Z/p\times \mu_p\times G$, where $G$ is local-local, and it follows
that $A$ is almost supersingular. Finally, if $\varphi$ is $(0,1,2)$, then the $p$-rank of $A$ is zero,
so the Newton polygon of $A$ is either supersingular or of $1/3$-type. \hfill $\Box$

\section{Open questions}

The main question of this paper is 

\bques
Let $A$ be an abelian variety of dimension $g$ over 
an algebraically closed field of characteristic $p>0$. What are the possible integers
$n_1\leqslant\ldots\leqslant n_r$ such that $\bU_A$ is isogenous to $\prod_{i=1}^r\bW_{n_i}$?
\eques

We can ask more specific questions about various particular types of abelian varieties.

The reader will notice that our proof of Theorem \ref{ss} in the non-supersingular case
bounded the dimension of $\End_\E(M^\vee/p^m,M/p^m)^\skew$ by the dimension of
the larger quasi-algebraic group $\End_\E(M^\vee/p^m,M/p^m)$. It is thus likely 
that in this case the bound $g-1$ for the $p$-exponent of $\bU_A$ is not the best possible.

\bques
Let $A$ be a non-supersingular abelian variety of dimension $g\geqslant 2$ over 
an algebraically closed field of characteristic $p>0$. Is it true that $p^{g-2}\bU_A=0$?
\eques

Assuming $p\neq 2$,
this holds for $g\leqslant 3$ (Proposition \ref{NPPabe3table}) and also for $g=4$ in the principally polarised case, as one sees by arguing similarly to Example \ref{d2}.

In the proof of Proposition \ref{3.6} we used a supersingular abelian variety $A$ with
Dieudonn\'e module $M=\E/(F^g-V^g)\E$. We note that the $a$-number of $A$
equals $\dim_k(M/(FM+VM))=1$, that is, $A$ is supergeneral.

\bques \label{Q2}
Let $A$ be a supersingular abelian variety of dimension $g\geqslant 3$ over 
an algebraically closed field of characteristic $p>0$ such that the $a$-number
of $A$ is at least $2$ (that is, $A$ is not supergeneral). Is it true that $p^{g-2}\bU_A=0$?
\eques

Here we must exclude the case $g=2$ because for a supersingular elliptic curve $E$ we have 
$a(E^2)=2$ whereas $\bU_{E^2}\cong\bG_{a,k}$.

Finally, one can enquire whether the description of the isogeny class of $\bU_A$ in Proposition \ref{3.6}
remains true for all supergeneral abelian varieties (this is \cite[Conjecture 5.2.1]{YuanThesis}).
In general it would be interesting to understand which numerical invariants of $M$ completely determine the isogeny type of $\bU_A$.

\bques \label{Q3}
Consider the abelian varieties $A$ over an algebraically closed field $k$ of characteristic $p>0$ 
with a fixed Newton polygon and $a$-number $1$. Is the isogeny type of $\bU_A$
 independent of $A$?
\eques

For $g\leqslant 3$
the answers to Questions \ref{Q2} and \ref{Q3} are positive, see Proposition \ref{NPPabe3table}.
The positive answer to Question \ref{Q3} for supergeneral varieties would be in agreement with Proposition \ref{supergeneralAp}, which says that in that case $\dim(\bU_A[p])=g-1$.

In the appendix below, we list the isogeny types of $\bU_A$
for some abelian varieties $A$ with $a$-number $1$. The following curious question seems
to have a positive answer for all the examples at our disposal.

\bques
Does the isogeny type $n_1\leqslant\ldots\leqslant n_r$ of $\bU_A$ always contain
a subsequence $1,2,\ldots,n_r-1,n_r$?
\eques

\appendix

\section{Numerical computations of isogeny types}

Below is a list of isogeny types of $\bU_A$ for a selected set of abelian varieties $A$. 

\def\arraystretch{1.5}
\begin{table}[h]
    \centering
\begin{tabular}{c|c c c}
$g$ & Newton polygon $\nu$   & $\dim(\bU_{A_\nu})$ & Isogeny type of $\bU_{A_\nu}$ \\ \hline
3 & $\left[\frac{1}{2}, \frac{1}{2}, \frac{1}{2}\right]$                                        & 3                          & $\left[1, 2\right]$          \\
  & $\left[\frac{1}{3}, \frac{2}{3}\right]$                                                     & 2                          & $\left[1, 1\right]$          \\ \hline
4 & $\left[\frac{1}{2}, \frac{1}{2}, \frac{1}{2}, \frac{1}{2}\right]$                           & 6                          & $\left[1, 2, 3\right]$       \\
  & $\left[\frac{1}{3}, \frac{1}{2}, \frac{2}{3}\right]$                                        & 4                          & $\left[1, 1, 2\right]$       \\
  & $\left[\frac{1}{4}, \frac{3}{4}\right]$                                                     & 3                          & $\left[1, 1, 1\right]$       \\ \hline
5 & $\left[\frac{1}{2}, \frac{1}{2}, \frac{1}{2}, \frac{1}{2}, \frac{1}{2}\right]$              & 10                         & $\left[1, 2, 3, 4\right]$    \\
  & $\left[\frac{2}{5}, \frac{3}{5}\right]$                                                     & 8                          & $\left[1, 2, 2, 3\right]$    \\
  & $\left[\frac{1}{3}, \frac{1}{2}, \frac{1}{2}, \frac{2}{3}\right]$                           & 7                          & $\left[1, 1, 2, 3\right]$    \\
  & $\left[\frac{1}{4}, \frac{1}{2}, \frac{3}{4}\right]$                                        & 5                          & $\left[1, 1, 1, 2\right]$    \\
  & $\left[\frac{1}{5}, \frac{4}{5}\right]$                                                     & 4                          & $\left[1, 1, 1, 1\right]$    \\ \hline
6 & $\left[\frac{1}{2}, \frac{1}{2}, \frac{1}{2}, \frac{1}{2}, \frac{1}{2}, \frac{1}{2}\right]$ & 15                         & $\left[1, 2, 3, 4, 5\right]$ \\
  & $\left[\frac{2}{5}, \frac{1}{2}, \frac{3}{5}\right]$                                        & 12                         & $\left[1, 2, 2, 3, 4\right]$ \\
  & $\left[\frac{1}{3}, \frac{1}{2}, \frac{1}{2}, \frac{1}{2}, \frac{2}{3}\right]$              & 11                         & $\left[1, 1, 2, 3, 4\right]$ \\
  & $\left[\frac{1}{3}, \frac{1}{3}, \frac{2}{3}, \frac{2}{3}\right]$                           & 10                         & $\left[1, 2, 2, 2, 3\right]$ \\
  & $\left[\frac{1}{4}, \frac{1}{2}, \frac{1}{2}, \frac{3}{4}\right]$                           & 8                          & $\left[1, 1, 1, 2, 3\right]$ \\
  & $\left[\frac{1}{5}, \frac{1}{2}, \frac{4}{5}\right]$                                        & 6                          & $\left[1, 1, 1, 1, 2\right]$ \\
  & $\left[\frac{1}{6}, \frac{5}{6}\right]$                                                     & 5                          & $\left[1, 1, 1, 1, 1\right]$
\end{tabular}
\end{table}
Here we consider abelian varieties of dimension $g$ with $p$-rank $0$ and $a$-number~$1$. Their
Dieudonné module is then isomorphic to  
\begin{equation}\label{forma}
    \E/(a_0F^g+a_1F^{g-1}+\dots+a_{2g}V^g)\E,
\end{equation}
for some $a_i\in W$, where $a_0,a_{2g}\in W^{\times}$. 

Let $\nu$ be a symmetric Newton polygon which has no slopes $0$ or $1$. If $M_{\nu}$ is the Dieudonné module of the form \eqref{forma} where
$$a_i=\begin{cases}
p^j, &(i,j) \text{ is a breaking point of }\nu\text{ and }i\le g;\\
p^{j+g-i}, &(i,j) \text{ is a breaking point of }\nu\text{ and }i> g;\\
0, &\text{ otherwise;}
\end{cases}$$ 
then $M_{\nu}$ has Newton polygon $\nu$ by \cite[Lemma 5.8]{LNV}. We take $A_{\nu}$ to be any abelian variety with Dieudonné module $M_{\nu}$.

The table contains the isogeny types of $\bU_{A_{\nu}}$ for $g\le 6$, when $p=3,5,7$ (the result 
does not depend on such choice of $p$). The computation was done in SageMath v10.7 using the implementation of Witt vectors in \cite{Witt}. The code can be found at \url{https://github.com/liviabasilico/brauer-p-torsion}. The principle of the SageMath script is to find the equations
that define $\Hom_\E(M^\vee/p^n,M/p^n)^\skew$ as a subvariety of 
$\Hom_{W_n}(M^\vee/p^n,M/p^n)$ and hence compute its dimension.

\noindent Institut de Recherche Mathématique Avancée, Université de Strasbourg, 7 rue René Descartes, 67000 Strasbourg, France

\texttt{livia.grammatica@math.unistra.fr}

\bigskip

\noindent Department of Mathematics, 
South Kensington Campus,
Imperial College London,
SW7~2AZ United Kingdom \  \ and \  \ 
Institute for the Information Transmission Problems,
Russian Academy of Sciences,
Moscow, 127994 Russia

\texttt{a.skorobogatov@imperial.ac.uk}

\bigskip

\noindent Beijing International Center for Mathematical Research,
Peking University, No. 5 Yiheyuan Road, Haidian District, Beijing, 100871 P.R. China

\texttt{yy}\verb|_|\texttt{@bicmr.pku.edu.cn}

\vspace{.5em}
\hrule
\vspace{.5em}
\begin{center}
\begin{minipage}{5cm}
\includegraphics[width=4.8cm]{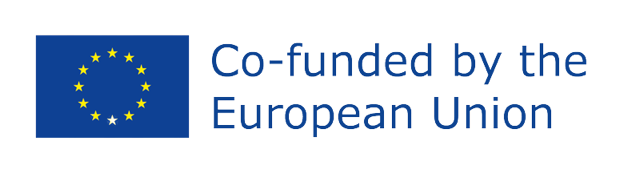}
\end{minipage}
\begin{minipage}{9.5cm}
\footnotesize Co-Funded by the European Union. Views and opinions expressed are however those of the author(s) only and do not necessarily reflect those of the European Union. Neither the European Union nor the granting authority can be held responsible for them.
\end{minipage}
\end{center}
\end{document}